\documentstyle[amsmath,amssymb,11pt]{article}

\renewcommand{\title}[1]{\leftline{\Large\bf #1}\par\medskip}
\renewcommand{\author}[1]{\medskip{\large #1}\par\medskip}

\newcommand{\sqplus}{\boxplus}

\newcommand{\rom}[1]{{\rm #1}}

\allowdisplaybreaks

\makeatletter\@addtoreset{equation}{section}\makeatother

\begin{document}

\setcounter{page}{1} \setcounter{section}{0} \thispagestyle{empty}

\newtheorem{definition}{Definition}[section]
\newtheorem{remark}{Remark}[section]
\newtheorem{proposition}{Proposition}[section]
\newtheorem{theorem}{Theorem}[section]
\newtheorem{corollary}{Corollary}[section]
\newtheorem{lemma}{Lemma}[section]

\newcommand{\D}{{\cal D}}
\newcommand{\N}{{\Bbb N}}
\newcommand{\C}{{\Bbb C}}
\newcommand{\Z}{{\Bbb Z}}
\newcommand{\R}{{\Bbb R}}
\newcommand{\Rp}{{\R_+}}
\newcommand{\eps}{\varepsilon}

\newcommand{\AS}{\operatorname{AS}}

\newcommand{\fii}{\varphi}

\newcommand{\supp}{\operatorname{supp}}
\newcommand{\la}{\langle}
\newcommand{\ra}{\rangle}
\newcommand{\const}{\operatorname{const}}
\newcommand{\ddGamma}{\overset{{.}{.}}{\Gamma}_X}

\newcommand{\ho}{\widehat\otimes}
\newcommand{\ot}{\otimes}

\renewcommand{\emptyset}{\varnothing}
\renewcommand{\tilde}{\widetilde}
\newcommand{\Formsg}{{\cal F}\Omega^{n+1}}
\newcommand{\Func}{{\cal FC}}

\newcommand{\di}{\partial}
\renewcommand{\div}{\operatorname{div}}

\begin{center}\LARGE \bf
LAPLACE OPERATORS\\[2mm]\LARGE\bf ON DIFFERENTIAL FORMS\\[2mm]
\LARGE\bf OVER CONFIGURATION SPACES\end{center}

\author{SERGIO ALBEVERIO}

\noindent{\sl Institut f\"{u}r Angewandte Mathematik,
Universit\"{a}t Bonn, Wegelerstr.~6, D-53115 Bonn, Germany;
\newline SFB 256, Univ.~Bonn, Germany;\newline SFB 237, Bochum--D\"usseldorf--Essen, Germany;
 \newline CERFIM (Locarno);
Acc.\ Arch.\ (USI), Switzerland;
\newline BiBoS, Univ.\ Bielefeld, Germany}\vspace{2mm}

\author{ALEXEI DALETSKII}

\noindent{\sl Institut f\"{u}r Angewandte Mathematik,
Universit\"{a}t Bonn, Wegelerstr.~6, D-53115 Bonn, Germany;
\newline SFB 256, Univ.~Bonn, Germany; \newline Institute of
Mathematics, Kiev, Ukraine;
\newline BiBoS, Univ.\ Bielefeld, Germany}\vspace{2mm}

\author{EUGENE LYTVYNOV}

\noindent{\sl Institut f\"{u}r Angewandte Mathematik,
Universit\"{a}t Bonn, Wegelerstr.~6, D-53115 Bonn, Germany;
\newline SFB 256, Univ.~Bonn, Germany;
\newline BiBoS, Univ.\ Bielefeld, Germany}

\begin{abstract}
Spaces of differential forms over configuration spaces with
Poisson measures are constructed. The corresponding Laplacians (of
Bochner and de Rham type) on forms and associated semigroups are
considered. Their probabilistic interpretation is given.
\end{abstract}

\noindent 2000 {\it AMS Mathematics Subject Classification}.
60G57, 58A10

\newpage
\tableofcontents

\section{Introduction}

Stochastic differential geometry of infinite-dimensional manifolds
has been a very active topic of research in recent times. One of
the important and intriguing problems discussed concerns the
construction of spaces of differential forms over such manifolds
and the study of the corresponding Laplace operators and
associated (stochastic) cohomologies. A central role in this
framework is played by the concept of the Dirichlet operator of a
differentiable measure, which is actually an infinite-dimensional
generalization of the Laplace--Beltrami operator on functions,
respectively the Laplace--Witten--de Rham operator on differential
forms. The study of the latter operator and the associated
semigroup on finite-dimensional manifolds was the subject of many
works, and it leads to deep results on the interface of stochastic
analysis, differential geometry and topology, and mathematical
physics, see, e.g.,\cite{E3}, \cite{EL}, \cite{CFKSi}, \cite
{Roe}, \cite{CoMo}. Dirichlet forms and processes in connection
with noncommutative $C^*$-algebras  were considered in e.g.\
\cite{Gr, AH-K, DL1}.

The interest in the infinite-dimensional case is motivated by
relations with supersymmetric quantum field theory. De Rham type
operators acting on differential forms over Hilbert spaces were
considered in \cite{Ar1}, \cite {Ar2}, \cite{ArM}, \cite{AK}. In
this relation, the mostly discussed example of an
infinite-dimensional non-flat space
  is the loop space of a
compact manifold, see \cite{W}, \cite{JL}, \cite{LRo}. Another
important example given by the infinite product of compact
manifolds was discussed in \cite{ADK1}, \cite{ADK2}, \cite{LeBe}.

At the same time, there is a growing interest in geometry and
analysis on Poisson spaces, i.e., on spaces of locally finite
configurations in noncompact manifolds equipped with the Poisson
measure. In \cite{AKR-1}, \cite {AKR0}, \cite{AKR1},  an approach
to these spaces as to infinite-dimensional manifolds  was
initiated. This approach is motivated by the connection of such
spaces with the theory of representations of diffeomorphism
groups, see \cite{GGPS}, \cite{VGG}, \cite {I} (these references
and \cite{AKR1}, \cite{AKR3} also contain discussion of relations
with quantum physics). In fact, the configuration space, which
does not possess  the structure of a smooth manifold in the proper
sense, can be equipped with some ``Riemannian-like'' structure
generated by the action of the diffeomorphism group of the initial
manifold. We refer the reader to \cite{AKR2}, \cite{AKR3},
\cite{Ro}, and references therein for further discussion of
analysis on Poisson spaces and applications.

In the present work, we develop this point of view. We define
spaces of differential forms over Poisson spaces and study Laplace
operators acting in these spaces. We show, in particular, that the
corresponding de Rham Laplacian can be expressed in terms of the
Dirichlet operator on functions on the Poisson space and the
Witten Laplacian on the initial manifold associated with the
intensity of the corresponding Poisson measure. We give a
probabilistic interpretation and investigate some properties of
the associated semigroups. The main general aim of our approach is
to develop a framework which extends to Poisson spaces (as
infinite-dimensional manifolds) the finite-dimensional Hodge--de
Rham theory.

The results of the present paper in the special case of 1-forms
were presented in \cite{ADL}.

A different approach to the construction of differential forms and
related objects over Poisson spaces, based on the ``transfer
principle'' from Wiener spaces, is proposed in \cite{Pr2}, see
also \cite{PPr} and \cite{Pr}.

\section{Differential forms over configuration spaces}

The aim of this section is to define differential forms over configuration
spaces (as infinite-dimensional manifold). First, we recall some known facts
and definitions concerning ``manifold-like'' structures and functional
calculus on these spaces.

\subsection{Functional calculus on configuration spaces}

Our presentation in this subsection is based upon \cite{AKR1},
however for later use in the present paper we give a different
description of some objects and results occurring in \cite{AKR1}.

Let $X$ be a complete, connected, oriented, $C^\infty $ (noncompact)
Riemannian manifold of dimension $d$. We denote by $\langle \bullet ,\bullet
\rangle _x$ the corresponding inner product in the tangent space $T_xX$ to $%
X $ at a point $x\in X$. The associated norm will be denoted by $|\bullet
|_x $ . Let also $\nabla ^X$ stand for the gradient on $X$.

The configuration space $\Gamma _X$ over $X$ is defined as the set of all
locally finite subsets (configurations) in $X$:
\begin{equation}
\Gamma _X:=\left\{ \,\gamma \subset X\mid |\gamma \cap \Lambda |<\infty
\text{ for each compact }\Lambda \subset X\,\right\} .  \nonumber
\end{equation}
Here, $|A|$ denotes the cardinality of the set $A$.

We can identify any $\gamma \in \Gamma _X$ with the positive integer-valued
Radon measure
\begin{equation}
\sum_{x\in \gamma }\varepsilon _x\in {\cal M}(X),  \nonumber
\end{equation}
where $\varepsilon _x$ is the Dirac measure with mass at $x$, $\sum_{x\in
\varnothing }\varepsilon _x:=$zero measure, and ${\cal M}(X)$ denotes the
set of all positive Radon measures on the Borel $\sigma $-algebra ${\cal B}
(X)$. The space $\Gamma _X$ is endowed with the relative topology as a
subset of the space ${\cal M}(X)$ with the vague topology, i.e., the weakest
topology on $\Gamma _X$ such that all maps
\begin{equation}
\Gamma _X\ni \gamma \mapsto \langle f,\gamma \rangle :=\int_Xf(x)\,\gamma
(dx)\equiv \sum_{x\in \gamma }f(x)  \nonumber
\end{equation}
are continuous. Here, $f\in C_0(X)$($:=$the set of all continuous functions
on $X$ with compact support). Let ${\cal B}(\Gamma _X)$ denote the
corresponding Borel $\sigma $-algebra.

Following \cite{AKR1}, we define the tangent space to $\Gamma _X$ at a point
$\gamma $ as the Hilbert space
\begin{equation}
T_\gamma \Gamma _X:=L^2(X\to TX;d\gamma ),  \nonumber
\end{equation}
or equivalently
\begin{equation}
T_\gamma \Gamma _X=\bigoplus_{x\in \gamma }T_xX  \label{tg-sp1}
\end{equation}
(compare also with \cite[Appendix 3]{VGG}).
The scalar product and the norm in $T_\gamma \Gamma _X$ will be denoted by $%
\langle \bullet ,\bullet \rangle _\gamma $ and $\left\| \bullet \right\|
_\gamma $, respectively. Thus, each $V(\gamma )\in T_\gamma \Gamma _X$ has
the form $V(\gamma )=(V(\gamma )_x)_{x\in \gamma }$, where $V(\gamma )_x\in
T_xX$, and
\begin{equation}
\| V(\gamma )\| _\gamma ^2=\sum_{x\in \gamma }|V(\gamma )_x|_x^2.  \nonumber
\end{equation}

Let $\gamma \in \Gamma _X$ and $x\in \gamma $. By ${\cal O}_{\gamma ,x}$ we
will denote an arbitrary open neighborhood of $x$ in $X$ such that the
intersection of the closure of ${\cal O}_{\gamma ,x}$ in $X$ with $\gamma
\setminus \{x\}$ is the empty set. For any fixed finite subconfiguration $%
\left\{ x_1,\dots ,x_k\right\} \subset \gamma $, we will always consider
open neighborhoods ${\cal O}_{\gamma ,x_1},\dots ,{\cal O}_{\gamma ,x_k}$
with disjoint closures.

Now, for a measurable function $F\colon\Gamma _X\to {\Bbb R}$, $\gamma \in
\Gamma _X$, and $\left\{ x_1,\dots ,x_k\right\} \subset \gamma $, we define
a function $F_{x_1,\dots ,x_k}(\gamma ,\bullet )\colon{\cal O}_{\gamma
,x_1}\times \dots \times {\cal O}_{\gamma ,x_k}\to {\Bbb R}$ by
\begin{multline*}
{\cal O}_{\gamma ,x_1}\times \dots \times {\cal O}_{\gamma ,x_k} \ni
(y_1,\dots ,y_k)\mapsto F_{x_1,\dots ,x_k}(\gamma ,y_1,\dots ,y_k):= \\
=F((\gamma \setminus \{x_1,\dots ,x_k\})\cup \{y_1,\dots ,y_k\})\in {\Bbb R}.
\end{multline*}
Since we will be interested only in the local behavior of the function $%
F_{x_1,\dots ,x_k}(\gamma ,\bullet )$ around the point $(x_1,\dots ,x_k)$,
we will not write explicitly which neighborhoods ${\cal O}_{\gamma ,x_i}$ we
use.

\begin{definition}
\label{def2.0}{\rm We say that a function $F:\Gamma _X\to {\Bbb R}$ is
differentiable at $\gamma \in \Gamma _X$ if for each $x\in \gamma $ the
function $F_x(\gamma ,\bullet )$ is differentiable at $x$ and
\[
\nabla ^\Gamma F(\gamma )=(\nabla ^\Gamma F(\gamma )_x)_{x\in \gamma }\in
T_\gamma \Gamma _X,
\]
where
\begin{equation}
\nabla ^\Gamma F(\gamma )_x:=\nabla ^XF_x(\gamma ,x).  \nonumber
\end{equation}
}
\end{definition}

{We will call $\nabla ^\Gamma F(\gamma )$ the {\it gradient\/} of
$F$ at $\gamma $. }

For a function $F$ differentiable at $\gamma $ and a vector
$V(\gamma )\in T_\gamma \Gamma _X$, the {\it directional
derivative\/} of $F$ at the point $\gamma $ along $V(\gamma )$ is
defined by
\begin{equation}
\nabla _V^\Gamma F(\gamma ):=\langle \nabla ^\Gamma F(\gamma ),V(\gamma
)\rangle _\gamma .  \nonumber
\end{equation}

In what follows, we will also use the shorthand notation
\begin{equation}
\nabla _x^XF(\gamma ):=\nabla ^XF_x(\gamma ,x),  \label{flick}
\end{equation}
so that
\[
\nabla ^\Gamma F(\gamma )=(\nabla _x^XF(\gamma ))_{x\in \gamma }.
\]
It is easy to see that the operation $\nabla ^\Gamma $ satisfies the usual
properties of differentiation, including the Leibniz rule.

We define a class ${\cal FC}$ of smooth cylinder functions on $\Gamma _X$ as
follows:

\begin{definition}
\label{def2.1} {\rm A measurable bounded function $F:\Gamma _X\to {\Bbb R}$
belongs to ${\cal FC}$ iff: }\newline
{\rm (i) there exists a compact $\Lambda \subset X$ such that $F(\gamma
)=F(\gamma _\Lambda )$ for all $\gamma \in \Gamma _X$, where $\gamma
_\Lambda :=\gamma \cap \Lambda $;}\newline
{\rm (ii) for any $\gamma \in \Gamma _X$ and $\left\{ x_1,\dots ,x_k\right\}
\subset \gamma $, $k\in {\Bbb N}$, the function $F_{x_1,\dots ,x_k}(\gamma
,\bullet )$ is infinitely differentiable with partial derivatives uniformly
bounded in $\gamma $ and $x_1,\dots ,x_k$ (i.e., the majorizing constant
depends only on the order of differentiation but not on the specific choice
of $\gamma \in \Gamma _X$, $k\in {\Bbb N}$, and $\{x_1,\dots ,x_k\}\subset
\gamma $).}
\end{definition}

Let us note that, for $F\in{\cal FC}$, only a finite number of coordinates
of $\nabla ^\Gamma F(\gamma )$ are not equal to zero, and so $\nabla ^\Gamma
F(\gamma )\in T_\gamma \Gamma _X$. Thus, each $F\in {\cal FC}$ is
differentiable at any point $\gamma \in \Gamma _X$ in the sense of
Definition~\ref{def2.0}.

\begin{remark}
\label{rem2.1}{\rm In \cite{AKR1}, the authors introduced the class ${\cal FC%
}_{\mathrm b}^\infty ({\cal D},\Gamma _X)$ of functions on $\Gamma
_X$ of the form
\begin{equation}
F(\gamma )=g_F(\left\langle \varphi _1,\gamma \right\rangle ,\dots
,\left\langle \varphi _N,\gamma \right\rangle ),  \label{2.1}
\end{equation}
where $g_F\in C_{\mathrm b}^\infty ({\Bbb R}^N)$ and $\varphi
_1,\dots ,\varphi _N\in {\cal D}:=C_0^\infty (X)$($:=$ the set of
all $C^\infty $-functions on $X$ with compact support). Evidently,
we have the inclusion
\[
{\cal FC}_{\mathrm b}^\infty ({\cal D},\Gamma _X)\subset {\cal
FC},
\]
and moreover, the gradient of $F$ of the form (\ref{2.1}) in the sense of
Definition~\ref{def2.0},
\begin{equation}
\nabla ^\Gamma F(\gamma )_x=\sum_{i=1}^N\frac{\partial g_F}{\partial s_i}%
(\langle \varphi _1,\gamma \rangle ,\dots ,\langle \varphi _N,\gamma \rangle
)\nabla ^X\varphi _i(x),  \nonumber
\end{equation}
coincides with the gradient of this function in the sense of \cite{AKR1}. }
\end{remark}

\subsection[Tensor bundles and cylinder forms over
 configuration spaces]{Tensor bundles and cylinder forms over configuration\\ spaces}

Our next aim is to introduce differential forms on $\Gamma _X$.

Vector fields and first order differential forms on $\Gamma _X$ will be
identified with sections of the bundle $T\Gamma _X.$ Higher order
differential forms will be identified with sections of tensor bundles $%
\wedge ^n(T\Gamma _X)$ with fibers
\begin{equation}
\wedge ^n(T_\gamma \Gamma _X)
\mbox{$:=$}%
\wedge ^n(L^2(X\rightarrow TX;\gamma )),  \nonumber
\end{equation}
where $\wedge ^n({\cal H})$ (or ${\cal H}^{\wedge n}$) stands for
the $n$-th antisymmetric tensor power of a Hilbert space ${\cal
H}$. In what follows, we will use different representations of
this space. Because of (\ref{tg-sp1}), we have
\begin{equation}
\wedge ^n(T_\gamma \Gamma _X)=\wedge ^n\bigg( \bigoplus_{x\in \gamma }T_xX %
\bigg) .  \label{tang-n}
\end{equation}

Let us introduce the factor space $X^n/S_n$, where $S_n$ is the permutation
group of $\{1,\dots ,n\}$ which naturally acts on $X^n$:
\begin{equation}
\sigma (x_1,\dots ,x_n)=(x_{\sigma (1)},\dots ,x_{\sigma (n)}),\qquad\sigma
\in S_n.  \nonumber
\end{equation}
The space $X^n/S_n$ consists of equivalence classes $[x_1,\dots ,x_n]$ and
we will denote by $[x_1,\dots ,x_n]_d$ an equivalence class $[x_1,\dots
,x_n] $ such that the equality $x_{i_1}=x_{i_2}=\dots =x_{i_k}$ can hold
only for $k\le d$ points. (In other words, any equivalence class $[x_1,\dots
,x_n]$ is a multiple configuration in $X$, while $[x_1,\dots ,x_n]_d$ is a
multiple configuration with multiplicity of points $\le d$.) In what
follows, instead of writing $[x_1,\dots ,x_n]_d:\{x_1,\dots ,x_n\}\subset
\gamma $, we will use the shortened notation $[x_1,\dots ,x_n]_d\subset
\gamma $, though $[x_1,\dots ,x_n]_d$ is not, of course, a set. We then have
from (\ref{tang-n}):
\begin{equation}
\wedge ^n(T_\gamma \Gamma _X)=\bigoplus_{[x_1,\dots
,x_n]_d\,\subset \gamma }T_{x_1}X\wedge T_{x_2}X\wedge \dots
\wedge T_{x_n}X.  \label{tang-n0}
\end{equation}
Here, the space $T_{x_1}X\wedge T_{x_2}X\wedge \dots \wedge
T_{x_n}X$ is understood as a subspace of the Hilbert space
$\big(T_{y_1}X\oplus T_{y_2}X\oplus\dots\oplus
T_{y_k}X\big)^{\otimes n},$ where $\{y_1,\dots,y_k\}$ is the set
of the different $x_j$'s, $j=1,\dots,n$. To see that
\eqref{tang-n0} holds, notice that
\begin{equation}\label{kjkjk}\big(T_{y_1}X\oplus
T_{y_2}X\oplus\dots\oplus T_{y_k}X\big)^{\otimes
n}\simeq\big(T_{y_{\nu(1)}}X\oplus
T_{y_{\nu(2)}}X\oplus\dots\oplus T_{y_{\nu(k)}}X\big)^{\otimes
n},\qquad \nu\in S_k\end{equation} (where $\simeq$ means
isomorphism), and moreover $T_{x_1}X\wedge
T_{x_2}X\wedge\dotsm\wedge T_{x_n}X$ and $T_{x_\sigma(1)}X\wedge
T_{x_\sigma(2)}X\wedge\dots\wedge T_{x_\sigma(n)}X$, $\sigma\in
S_n$, coincide   as  subspaces of the space \eqref{kjkjk}.

Thus, under a differential form $W$ of order $n$, $n\in {\Bbb N}$, over $%
\Gamma _X,$ we will understand a mapping
\begin{equation}  \label{lkghf}
\Gamma _X\ni \gamma \mapsto W(\gamma )\in \wedge ^n(T_\gamma \Gamma _X).
\end{equation}
We denote by $W (\gamma )_{[x_1,\dots ,x_n]_d}$ the corresponding component
of $W(\gamma)$ in the decomposition (\ref{tang-n0}).

In particular, in the case $n=1$, a 1-form $V$ over $\Gamma _X$ is given by
the mapping
\begin{equation}
\Gamma _X\ni \gamma \mapsto V(\gamma )=(V(\gamma )_x)_{x\in \gamma }\in
T_\gamma \Gamma _X.  \nonumber
\end{equation}

For fixed $\gamma \in \Gamma _X$ and $x\in \gamma $, we consider the mapping
\begin{equation}
{\cal O}_{\gamma ,x}\ni y\mapsto W_x(\gamma ,y)%
\mbox{$:=$}%
W(\gamma _y)\in \wedge ^n(T_{\gamma _y}\Gamma _X),  \nonumber
\end{equation}
where $\gamma _y%
\mbox{$:=$}%
(\gamma \setminus \{x\})\cup \{y\},$ which is a section of the Hilbert
bundle
\begin{equation}
\wedge ^n(T_{\gamma _y}\Gamma _X)\mapsto y\in {\cal O}_{\gamma ,x}
\label{bund1}
\end{equation}
over ${\cal O}_{\gamma ,x}.$ The Levi--Civita connection on $TX$ generates
in a natural way a ``product'' connection on this bundle. We denote by $%
\nabla _{\gamma ,x}^X$ the corresponding covariant derivative, and use the
notation
\begin{equation}
\nabla _x^XW(\gamma )%
\mbox{$:=$}%
\nabla _{\gamma ,x}^X\,W_x(\gamma ,x)\in T_xX\otimes \left( \wedge
^n(T_\gamma \Gamma _X)\right)  \nonumber
\end{equation}
if the section $W_x(\gamma ,\bullet )$ is differentiable at $x$.
Analogously, we denote by $\Delta _x^X$ the corresponding Bochner Laplacian
associated with the volume measure $m$ on ${\cal O}_{\gamma ,x}$ (see
subsec.~3.2 where the notion of Bochner Laplacian is recalled).

Similarly, for a fixed $\gamma \in \Gamma _X$ and $\left\{
x_1,\dots,x_k\right\} \subset \gamma $, we define a mapping
\begin{multline*}
{\cal O}_{\gamma ,x_1}\times \dots\times {\cal O}_{\gamma ,x_k}
\ni (y_1,\dots,y_k)\mapsto W _{x_1,\dots,x_k}(\gamma
,y_1,\dots,y_k) := \\ = W(\gamma _{y_1,\dots,y_k}) \in \wedge
^n(T_{\gamma _{y_1,\dots,y_k}}\Gamma _X),
\end{multline*}
where $\gamma _{y_1,\dots,y_k}
\mbox{$:=$}%
(\gamma \setminus \{x_1,\dots,x_k\})\cup \{y_1,\dots,y_k\}$, which
is a section of the Hilbert bundle
\begin{equation}
\wedge ^n(T_{\gamma _{y_1,\dots,y_k}}\Gamma _X)\mapsto \left(
y_1,\dots,y_k\right) \in {\cal O}_{\gamma ,x_1}\times \dots\times
{\cal O} _{\gamma ,x_k}  \label{bund-n}
\end{equation}
over ${\cal O}_{\gamma ,x_1}\times\dots\times {\cal O}_{\gamma
,x_k}.$

Let us remark that, for any $\eta \subset \gamma $, the space $\wedge
^n(T_\eta \Gamma _X)$ can be identified in a natural way with a subspace of $%
\wedge ^n(T_\gamma \Gamma _X)$. In this sense, we will use expressions of
the type $W (\gamma )=W (\eta )$ without additional explanations.

A set ${\cal F}\Omega ^n$ of smooth cylinder $n$-forms over $\Gamma _X$ will
be defined as follows.

\begin{definition}
\label{def2.2}{\rm ${\cal F}\Omega ^n$ is the set of $n$-forms $W$ over $%
\Gamma _X$ which satisfy the following conditions: }\newline
{\rm (i) there exists a compact $\Lambda =\Lambda (W)\subset X$ such that $%
W(\gamma )=W(\gamma _\Lambda )$; }\newline {\rm (ii) for each
$\gamma \in \Gamma _X$ and $\left\{ x_1,\dots ,x_k\right\} \subset
\gamma $, the section $W_{x_1,\dots ,x_k}(\gamma
,\bullet )$ of the bundle (\ref{bund-n}) is infinitely differentiable at $%
(x_1,\dots ,x_k)$, and bounded together with partial derivatives
component-wise in the sense of decomposition \eqref{tang-n0},
uniformly in $\gamma , x_1,\dots,x_k$, and the component. }
\end{definition}

\begin{remark}
\label{form-fin}{\rm For each $W\in {\cal F}\Omega ^n$, $\gamma \in \Gamma
_X $, and any open bounded $\Lambda \supset \Lambda (W)$, we can define the
form $W_{\Lambda ,\gamma }$ on ${\cal O}_{\gamma ,x_1}\times \dots \times
{\cal O}_{\gamma ,x_k}$ by
\begin{equation}
W_{\Lambda ,\gamma }(y_1,\dots ,y_k)=\operatorname{Proj}_{\wedge
^n(T_{y_1}X\oplus \dots \oplus T_{y_k}X)}W((\gamma \setminus \{x_1,\dots
,x_k\})\cup \{y_1,\dots ,y_k\}),  \label{cyl-form}
\end{equation}
where $\{x_1,\dots ,x_k\}=\gamma \cap \Lambda $. The item (ii) of
Definition~ \ref{def2.2} is obviously equivalent to the assumption that $%
W_{\Lambda ,\gamma }$ is smooth and bounded together with all
partial derivatives component-wise uniformly in $\gamma $ (for
some $\Lambda $ and consequently for any $\Lambda \supset \Lambda
(W)$).}
\end{remark}

\begin{definition}
\label{def2.3}{\rm We define the covariant derivative $\nabla ^\Gamma W$ of
a form $W$ given by \eqref{lkghf} as the mapping
\begin{equation}
\Gamma _X\ni \gamma \mapsto \nabla ^\Gamma W(\gamma )%
\mbox{$:=$}%
(\nabla _x^XW(\gamma ))_{x\in \gamma }\in T_\gamma \Gamma _X\otimes \left(
\wedge ^n(T_\gamma \Gamma _X)\right)  \nonumber
\end{equation}
if for all $\gamma \in \Gamma _X$ and $x\in \gamma $ the form $W_x(\gamma
,\bullet )$ is differentiable at $x$ and the $\nabla ^\Gamma W(\gamma )$
just defined indeed belongs to $T_\gamma \Gamma _X\otimes \left( \wedge
^n(T_\gamma \Gamma _X)\right) $. }
\end{definition}

\begin{remark}
{\rm For each $W\in {\cal F\Omega }^n$, the covariant derivative $\nabla
^\Gamma W$ exists, and moreover only a finite number of the coordinates $%
\nabla ^\Gamma W(\gamma )_{x,\,[x_1,\dots ,x_n]_d}$ in the decomposition
\[
T_\gamma \Gamma _X\otimes \big( \wedge ^n(T_\gamma \Gamma _X)\big)%
=\bigoplus_{x\in \gamma ,\,[x_1,\dots ,x_n]_d\subset \gamma }T_xX\otimes
(T_{x_1}X\wedge \dots \wedge T_{x_n}X)
\]
are not equal to zero.}
\end{remark}

\begin{remark}
\label{vbnk}{\rm For each $W\in {\cal F}\Omega ^n$, $\gamma \in \Gamma _X$, $%
x\in \gamma $, and $[x_1,\dots ,x_n]_d\subset \gamma $, we define the
mapping $W_x(\gamma ,\bullet )_{[x_1,\dots ,x_n]_d}$ as follows: if $x\ne
x_j $ for all $j=1,\dots ,n$, then
\[
{\cal O}_{\gamma ,x}\ni y\mapsto W_x(\gamma ,y)_{[x_1,\dots
,x_n]_d}:=W((\gamma \setminus \{x\})\cup \{y\})_{[x_1,\dots ,x_n]_d}\in
T_{x_1}X\wedge \dots \wedge T_{x_n}X,
\]
and if $x=x_i$ for some $x_i\in \{x_1,\dots ,x_n\}$, then
\[
{\cal O}_{\gamma ,x}\ni y\mapsto W_x(\gamma ,y)_{[x_1,\dots
,x_n]_d}:=W((\gamma \setminus \{x\})\cup \{y\})_{[y_1,\dots ,y_n]_d}\in
T_{y_1}X\wedge \dots \wedge T_{y_n}X,
\]
where $y_j=x_j$ if $x\ne x_j$ and $y_j=y$ otherwise. Then, the condition
(ii) of Definition~\ref{def2.2} yields, in particular, that the mapping $%
W_x(\gamma ,\bullet )_{[x_1,\dots ,x_n]_d}$ is $C^\infty $ for all $x\in
\gamma $ and $[x_1,\dots ,x_n]_d\subset \gamma $. Now, we have
\[
\nabla ^\Gamma W(\gamma)_{x,\,[x_1,\dots ,x_n]_d}=\nabla
_x^XW(\gamma )_{[x_1,\dots ,x_n]_d},
\]
where
\[
\nabla _x^XW(\gamma )_{[x_1,\dots ,x_n]_d}:=\nabla ^XW_x(\gamma
,x)_{[x_1,\dots ,x_n]_d}.
\]
Notice that, in the case where $x\ne x_j$ for all $j=1,\dots ,n$, $\nabla
^XW_x(\gamma ,\bullet )_{[x_1,\dots ,x_n]_d}$ means, in fact, the usual
derivative of a mapping defined on ${\cal O}_{\gamma ,x}$ and taking values
in the fixed vector space $T_{x_1}X\wedge \dots \wedge T_{x_n}X$. On the
other hand, if $x$ does coincide with some $x_i\in \{x_1,\dots ,x_n\}$, then
the expression $\nabla ^XW_x(\gamma ,x)_{[x_1,\dots ,x_n]_d}$ can be
understood as the $T_xX\otimes (T_{x_1}X\wedge \dots \wedge T_{x_n}X)$%
-coordinate of the covariant derivative of the $n$-form
\begin{multline}
{\cal O}_{\gamma ,y_1}\times \dots \times {\cal O}_{\gamma ,y_k}\ni
(z_1,\dots ,z_k)\mapsto  \label{2.4} \\
\mapsto \operatorname{Proj}_{\wedge ^n(T_{z_1}X\oplus \dots \oplus
T_{z_k}X)}W((\gamma \setminus \{y_1,\dots ,y_k\})\cup \{z_1,\dots ,z_k\})
\end{multline}
at the point $(y_1,\dots y_k)$, where $\{y_1,\dots ,y_k\}$ is the set of all
the different $x_j$'s, $j=1,\dots ,n$. In fact, the last sentence was just
an alternative description of the notion of covariant derivative $\nabla
^XW_x(\gamma ,x)_{[x_1,\dots ,x_n]_d}$ of the mapping $W_x(\gamma ,\bullet
)_{[x_1,\dots ,x_n]_d}$ in the case where $x$ coincides with some $x_i$. }
\end{remark}

\begin{proposition}
\label{prop2.1} For arbitrary $W^{(1)},W^{(2)}\in {\cal F}\Omega ^n$, we
have
\begin{multline*}
\nabla ^\Gamma \langle W^{(1)}(\gamma ),W^{(2)}(\gamma )\rangle _{\wedge
^n(T_\gamma \Gamma _X)}= \\
=\langle \nabla ^\Gamma W^{(1)}(\gamma ),W^{(2)}(\gamma )\rangle _{\wedge
^n(T_\gamma \Gamma _X)}+\langle W^{(1)}(\gamma ),\nabla ^\Gamma
W^{(2)}(\gamma )\rangle _{\wedge ^n(T_\gamma \Gamma _X)}.
\end{multline*}
\end{proposition}

\noindent
{\it Proof}. We have, for any fixed $\gamma \in \Gamma _X$,
\begin{gather*}
\nabla ^\Gamma \langle W ^{(1)}(\gamma ),W ^{(2)}(\gamma )\rangle _{\wedge
^n(T_\gamma \Gamma _X)} =\sum_{x\in\gamma} \nabla^X_x\langle W
^{(1)}(\gamma),W^{(2)}(\gamma)\rangle_{\wedge^n(T_\gamma\Gamma_X)} \\
=\sum_{x\in\gamma}\nabla^X_x\sum_{[x_1,\dots,x_n]_d\subset\gamma}\langle
W^{(1)}(\gamma)_{[x_1,\dots,x_n]_d},W^{(2)}(\gamma)_{[x_1,\dots,x_n]_d}
\rangle_{T_{x_1}X\wedge\dots\wedge T_{x_n}X} \\
=\sum_{x\in\gamma}\,\sum_{[x_1,\dots,x_n]_d\subset\gamma}\nabla^X_x\langle
W^{(1)}(\gamma)_{[x_1,\dots,x_n]_d},W^{(2)}(\gamma)_{[x_1,\dots,x_n]_d}
\rangle_{T_{x_1}X\wedge\dots\wedge T_{x_n}X} \\
=\sum_{x\in\gamma}\,\sum_{[x_1,\dots,x_n]_d\subset\gamma} \big[ %
\langle\nabla^X_xW^{(1)}(\gamma)
_{[x_1,\dots,x_n]_d},W^{(2)}_{[x_1,\dots,x_n]_d}\rangle_{T_{x_1}X\wedge
\dots\wedge T_{x_n}X} \\
\text{}+\langle W^{(1)}(\gamma)_{[x_1,\dots,x_n]_d}
,\nabla^X_xW^{(2)}(\gamma)_{[x_1,\dots,x_n]_d}\rangle_{T_{x_1}X\wedge\dots
\wedge T_{x_n}X}\big] \\
=\langle\nabla ^\Gamma W ^{(1)}(\gamma ),W ^{(2)}(\gamma )\rangle _{\wedge
^n(T_\gamma \Gamma _X)}+\langle W ^{(1)}(\gamma ),\nabla ^\Gamma
W^{(2)}(\gamma )\rangle _{\wedge ^n(T_\gamma \Gamma _X)}.
\end{gather*}
(All the sums above are actually finite because of the definition of ${\cal %
F }\Omega^n$.)\quad$\blacksquare$

\subsection{Square-integrable $n$-forms}

Our next goal is to give a description of the space of $n$-forms that are
square-integrable with respect to the Poisson measure.

Let $m$ be the volume measure on $X$, let $\rho \colon X\to {\Bbb
{R}}$ be a measurable function such that $\rho >0$ $m$-a.e., and
$\rho ^{1/2}\in H_{\mathrm loc}^{1,2}(X)$, and define the measure
$\sigma (dx):=\rho (x)\,m(dx)$.
Here, $H_{\mathrm loc}^{1,2}(X)$ denotes the local Sobolev space of order 1 in $%
L_{\mathrm loc}^2(X;m)$. Then, $\sigma $ is a nonatomic Radon
measure on $X$.

Let $\pi _\sigma $ stand for the Poisson measure on $\Gamma _X$ with
intensity $\sigma $. This measure is characterized by its Laplace transform
\begin{equation}
\int_{\Gamma _X}e^{\langle f,\gamma \rangle }\,\pi _\sigma (d\gamma )=\exp
\int_X(e^{f(x)}-1)\,\sigma (dx),\qquad f\in {\cal D}.  \nonumber
\end{equation}
Let $F\in L^1(\Gamma _X;\pi _\sigma )$ be cylindrical, that is, there exits
a compact $\Lambda \subset X$ such that $F(\gamma )=F(\gamma _\Lambda )$.
Then, one has the following formula, which we will use many times:
\begin{equation}
\int_{\Gamma _X}F(\gamma )\,\pi _\sigma (d\gamma )=e^{-\sigma (\Lambda
)}\sum_{n=0}^\infty \frac 1{n!}\int_{\Lambda ^n}F(\{x_1,\dots
,x_n\})\,\sigma (dx_1)\dotsm\sigma (dx_n).  \label{3.1}
\end{equation}

We define on the set ${\cal F}\Omega^n$ the $L^2$-scalar product with
respect to the Poisson measure:
\begin{equation}
(W^{(1)},W^{(2)})_{L_{\pi _\sigma }^2\Omega^n}%
\mbox{$:=$}%
\int_{\Gamma _X}\langle W ^{(1)}(\gamma ),W ^{(2)}(\gamma )\rangle _{\wedge
^nT_\gamma \Gamma _X}\,\pi _\sigma (d\gamma ).  \label{4.1}
\end{equation}
As easily seen, for each $W \in {\cal F}\Omega ^n$, there exists $\varphi\in%
{\cal D}$, $\varphi\ge0$, such that
\[
|\langle W(\gamma ),W(\gamma )\rangle _{\wedge ^nT_\gamma \Gamma
_X}|\le\langle \varphi^{\otimes n},\gamma^{\otimes n}\rangle.
\]
Hence, the function under the sign of integral in \eqref{4.1} indeed belongs
to $L^1(\Gamma _X;\pi _\sigma )$, since the Poisson measure has all moments
finite. Moreover, $(W ,W )_{L_{\pi _\sigma }^2\Omega ^n}>0$ if $W $ is not
identically zero. Hence, we can define the Hilbert space
\begin{equation}
L_{\pi _\sigma }^2\Omega ^n%
\mbox{$:=$}%
L^2(\Gamma _X\to \wedge ^n(T\Gamma _X);\pi _\sigma )  \nonumber
\end{equation}
as the closure of ${\cal F}\Omega ^n$ in the norm generated by the scalar
product (\ref{4.1}).

We will give now an isomorphic description of the space $L_{\pi _\sigma
}^2\Omega ^n$ via the space $L^2_{\pi_\sigma}(\Gamma_X):= L^2(\Gamma_X\to
{\Bbb {R};\pi_\sigma)}$ and some special spaces of square-integrable forms
on $X^m$, $m=1,\dots,n$.


We need first some preparations. Let $X^m$ be the $m$-th Cartesian power of
the manifold $X$. We have
\begin{equation}
\wedge ^n(T_{(x_1,\dots ,x_m)}X^m)=\bigoplus_{\begin{gathered}{\scriptstyle{
0\le k_1,\dots,k_m\le d}} \\ \scriptstyle k_1+\dots+k_m=n \end{gathered}%
}(T_{x_1}X)^{\wedge k_1}\wedge \dots \wedge (T_{x_m}X)^{\wedge k_m}.
\label{n-forms}
\end{equation}
For an $n$-form $\omega $ on $X^m$, we denote by $\omega (x_1,\dots
,x_m)_{k_1,\dots ,k_m}$ the corresponding component of $\omega (x_1,\dots
,x_m)$ in the decomposition (\ref{n-forms}).

Let {\rm \
\[
\widetilde{X}^m:=\big\{\,\bar{x}=(x_1,\dots,x_m)\in X^m\mid x_i\ne
x_j\text{ if }i\ne j\,\big\}.
\]
}We introduce a set $\Psi _0^n(\widetilde{X}^m)$ (resp.\ $\Psi
_0^n(X^m)$) of bounded $n$-forms $\omega $ over $X^m$ which have
compact support, smooth on $\widetilde{X}^m$ (resp.\ on $X^m$),
and satisfy the following assumptions:

\begin{enumerate}
\item[(i)]  $\omega (x_1,\dots ,x_m)_{k_1,\dots ,k_m}=0$ if $k_j=0$ for some
$j\in \{1,\dots ,m\}$;

\item[(ii)]  $\omega $ is symmetric:
\begin{equation}
\omega (x_1,\dots ,x_m)=\omega (x_{\sigma (1)},\dots ,x_{\sigma (m)})\qquad
\text{for each }\sigma \in S_m.  \label{symmetric}
\end{equation}
(we identify the spaces $\wedge ^n(T_{(x_1,\dots ,x_m)}X^m)$ and
$\wedge ^n(T_{(x_{\sigma (1)},\dots ,x_{\sigma (m)})}X^m)$, see
(\ref{n-forms}) and the explanation just after formula
\eqref{tang-n0}).
\end{enumerate}

Evidently, $\Psi _0^n(X^m)\subset \Psi _0^n(\widetilde{X}^m)$.

Let ${:}\,\gamma ^{\otimes m}\,{:}$ be the measure on $X^m$ given by
\begin{equation}
{:}\,\gamma ^{\otimes m}\,{:}(dx_1,\dots ,dx_m):=\sum_{\left\{ y_1,\dots
,y_m\right\} \subset \gamma }\varepsilon _{y_1}\widehat{\otimes }\dotsm%
\widehat{\otimes }\varepsilon _{y_m}(dx_1,\dots ,dx_m),  \nonumber
\end{equation}
where
\[
\varepsilon _{y_1}\widehat{\otimes }\dotsm\widehat{\otimes }\varepsilon
_{y_m}(dx_1,\dots ,dx_m):=\frac 1{m!}\sum_{\sigma \in s_m}\varepsilon
_{y_{\sigma (1)}}\otimes \dots \otimes \varepsilon _{y_{\sigma
(m)}}(dx_1,\dots ,dx_m).
\]
 We will use the notation
\begin{equation}
{\Bbb T}_{\{x_1,\dots ,x_m\}}^{(n)}X^m:=\bigoplus_{\begin{gathered}{%
\scriptstyle{ 1\le k_1,\dots,k_m\le d}} \\ \scriptstyle k_1+\dots+k_m=n
\end{gathered}}(T_{x_1}X)^{\wedge k_1}\wedge \dots \wedge (T_{x_m}X)^{\wedge
k_m}.  \label{tang-n2}
\end{equation}
By virtue of \eqref{tang-n0}, we have
\begin{equation}
\wedge ^n(T_\gamma \Gamma
_X)=\bigoplus_{m=1}^n\bigoplus_{\{x_1,\dots ,x_m\}\subset \gamma
}{\Bbb T}_{\{x_1,\dots ,x_m\}}^{(n)}X^m.  \label{tang-n1}
\end{equation}
For $W\in {\cal F}\Omega ^n$, we denote by $W_m(\gamma )\in
\bigoplus_{\{x_1,\dots ,x_m\}\subset \gamma }{\Bbb T}_{\{x_1,\dots
,x_m\}}^{(
n)}X^m$ the corresponding component of $W(\gamma )$ in the decomposition %
\eqref{tang-n1}. Thus, for $\{x_1,\dots ,x_m\}\subset \gamma $,
$W_m(\gamma )(x_1,\dots ,x_m)$ is equal to the projection of
$W(\gamma )\in \wedge ^n(T_\gamma \Gamma _X)$ onto the subspace
${\Bbb T}_{\{x_1,\dots ,x_m\}}^{(n)}X^m$.

For $\bar{x}=(x_1,\dots,x_m)\in \widetilde{X}^m$  we set  $\{\bar
x\}:=\{x_1,\dots,x_m\}$.

\begin{lemma}
For $W,V\in {\cal F}\Omega ^n$, we have
\begin{equation}
\left( W\left( \gamma \right) ,V\left( \gamma \right) \right)
_{\wedge ^n(T_\gamma \Gamma
_X)}=\sum_{m=1,\dots,n}\int_{X^m}\left( W_m\left( \gamma \right)
\left( \bar{x}\right) ,V_m\left( \gamma \right) \left(
\bar{x}\right) \right) _{{\Bbb T}_{\left\{ \bar{x}\right\}
}^{(n)}X^m}{:}\,\gamma ^{\otimes m}\,{:}\left( d\bar{x}\right)
\label{tensor-ng}
\end{equation}
\end{lemma}

The proof can be obtained by a direct calculation.

Let us remark that each $\omega \in \Psi _0^n(\widetilde{X}^m)$ generates a
cylinder form $W\in {\cal F}\Omega ^n$ by the formula
\begin{equation}
W_k(\gamma )(x_1,\dots ,x_k)=\begin{cases}\omega(x_1,\dots,x_m),&k=m,\\
0,&k\ne m.\end{cases}  \nonumber  \label{const-form}
\end{equation}

Let us denote by $L_\sigma ^2\Psi _0^n(X^m)$ the space obtained as the
completion of $\Psi _0^n(X^m)$ in the $L^2$-scalar product w.r.t.\ the
measure $\sigma ^{\otimes m}$. Evidently, $\Psi _0^n(\widetilde{X}^m)$ is a
dense subset of $L_\sigma ^2\Psi _0^n(X^m).$ We have

\begin{proposition}
\label{sq-int}The space $L_{\pi _\sigma }^2\Omega ^n$ is unitarily
isomorphic to the space
\begin{equation}
L_{\pi _\sigma }^2(\Gamma _X)\otimes \bigg[ \bigoplus_{m=1}^nL_\sigma ^2\Psi
^n(X^m)\bigg]=\bigoplus_{m=1}^nL_{\pi _\sigma }^2(\Gamma _X)\otimes L_\sigma
^2\Psi ^n(X^m),  \label{tensor-n}
\end{equation}
where the corresponding isomorphism $I^n$ is defined by the formula
\begin{equation}
I_m^nW(\gamma ,\bar{x})%
\mbox{$:=$}%
(m!)^{-1/2}\,W_m(\gamma \cup \{\bar{x}\})(\bar{x}),\qquad
m=1,\dots ,n. \label{asfdgfsdd}
\end{equation}
Here{\rm ,} $I_m^nW:=\left( I^nW\right) _m$ is the $m$-th component of $%
I^nV $ in the decomposition {\rm (\ref{tensor-n}).}
\end{proposition}

\begin{remark}
{\rm Actually, the formula \eqref{asfdgfsdd} makes sense only for $\bar{%
x}\in \widetilde{X}^m$. However, since the set $X^m\setminus
\widetilde{X}^m$ is of zero $\sigma ^{\otimes m}$ measure, this
does not lead to a contradiction.}
\end{remark}

\noindent{\it Proof}. First, we recall an extension of the Mecke
identity (e.g.\ \cite{KMM}) to the case of functions of several
variables \cite{Gena}:
\begin{multline}
\int_{\Gamma _X}\bigg[ \int_{X^m}f(\gamma ,\bar{x})\,{:}\,\gamma
^{\otimes m}\,{:}\,(d\bar{x})\bigg] \pi _\sigma (d\gamma ) \\
=\frac1{m!}\,\int_{\Gamma _X}\bigg[ \int_{X^m}f(\gamma \cup \{\bar{x}\},\bar{x}%
)\,\sigma ^{\otimes m}(d\bar{x})\bigg] \pi _\sigma (d\gamma ),
\label{mecke-n}
\end{multline}
where $f\colon\Gamma _X\times X^m\to {\Bbb {R}}$ is a measurable
function for which at least one of the double-integrals in
\eqref{mecke-n} exists (this formula can be easily proved by a
direct calculation using \eqref{3.1} for $f(\gamma
,\bar{x})=F(\gamma )g(\bar{x})$, where $F(\gamma )$ is bounded and
cylindrical and $g(\bar{x})$ is bounded and has compact support).

Next, let us specify the scalar product of two cylinder $n$-forms $W,V\in
{\cal F}\Omega ^n.$ We have, according to (\ref{tensor-ng}),
\begin{align}
& \left( W(\gamma ),V(\gamma )\right) _{\wedge ^n(T_\gamma \Gamma
_X)}= \nonumber \\ & \qquad =\sum_{m=1}^n\int_{X^m}\left(
W_m(\gamma )(\bar{x}),V_m(\gamma
)(\bar{x})\right) _{\bar{x}}\,{:}\,\gamma ^{\otimes m}\,{:}\,(d%
\bar{x})  \nonumber \\
& \qquad =\sum_{m=1}^n\int_{X^m}(W_m(\gamma \cup \{\bar{x}\})(\bar{%
x}),V_m(\gamma \cup \{\bar{x}\})(\bar{x}))_{\bar{x}}\,{:}%
\,\gamma ^{\otimes m}\,{:}\,(d\bar{x}),
\end{align}
where $\left( \bullet ,\bullet \right) _{{\bar{x}}}%
\mbox{$:=$}%
\left( \bullet ,\bullet \right) _{{\Bbb T}_{\left\{
\bar{x}\right\} }^{(n)}X^m}$ (we used the evident equality $\gamma
\cup \{\bar{x}\}=\gamma $ for $\{\bar{x}\}\subset \gamma $). The
application of the Mecke identity (\ref{mecke-n}) to the function
\[
f(\gamma ,\bar{x})=\left( W_m(\gamma \cup \{\bar{x}\})(\bar{x}%
),V_m(\gamma \cup \{\bar{x}\})(\bar{x})\right) _{\bar{x}}
\]
shows that
\begin{equation}
\left( W,V\right) _{L_{\pi _\sigma }^2\Omega ^n}=
\sum_{{m=1}}^n  \frac1{m!}\,   \int_{\Gamma _X}\int_{X^m}\left( W_m(\gamma \cup \{\bar{x%
}\})(\bar{x}),V_m(\gamma \cup \{\bar{x}\})(\bar{x})\right) _{%
\bar{x}}\,\sigma ^{\otimes m}(d\bar{x} )\,\pi _\sigma (d\gamma ).
\nonumber
\end{equation}
The space ${\cal F}\Omega ^n$ is dense in $L_{\pi _\sigma }^2\Omega ^n$, and
so it remains only to show that $I^n({\cal F}\Omega ^n)$ is a dense subspace
of $\bigoplus_{m=1}^nL_{\pi _\sigma }^2(\Gamma _X)\otimes L_\sigma ^2\Psi
^n(X^m)$, i.e., $I_m^n({\cal F}\Omega ^n)$ is a dense subspace of $L_{\pi
_\sigma }^2(\Gamma _X)\otimes L_\sigma ^2\Psi ^n(X^m)$, $m=1,\dots ,n$.

For $F\in {\cal FC}$ and $\omega \in \Psi _0^n({X}^m)$, we define a form $W$
by setting
\begin{equation}
\begin{gathered}W_k(\gamma):=0\qquad\text{for }k\ne m,\\
W_m(\gamma)(\bar x):= (m!)^{1/2}\,F(\gamma\setminus\{\bar
x\})\omega(\bar x).\end{gathered}  \label{kikimora}
\end{equation}
Evidently, we have $W\in {\ {\cal F}}\Omega ^n$ and
\begin{equation}
\begin{gathered}I^n_kW(\gamma,\bar x)=0\qquad\text{for }k\ne m,\\
I_m^nW(\gamma ,\bar{x})=F(\gamma )\omega (\bar{x})\end{gathered}
\label{is1}
\end{equation}
for each $\gamma \in \Gamma _X$ and $\bar{x}\in \tilde{X}^m$ such that $%
\{\bar{x}\}\cap \gamma =\varnothing $. Since $\gamma $ is a set of
zero $\sigma $ measure, we conclude from (\ref{is1}) that
\[
I_m^nW=F\otimes \omega .
\]
Noting that the linear span of such $F\otimes \omega $ is dense in $L_{\pi
_\sigma }^2(\Gamma _X)\otimes L_\sigma ^2\Psi ^n(X^m)$, we obtain the
result. \vspace{2mm}\quad ${\blacksquare }$

In what follows, we will denote by ${\cal D}\Omega ^n$ the linear span of
forms $W$ defined by (\ref{kikimora}), $m=1,\dots,n$. As we already noticed
in the proof of Proposition \ref{sq-int}, ${\cal D}\Omega ^n$ is a subset of
${\cal F} \Omega ^n$ and is dense in $L_{\pi _{\text{\/}\sigma }}^2\Omega ^n$%
.

\begin{corollary}
\label{fock} We have the unitary isomorphism
\[
{\cal I}^n\colon L_{\pi _\sigma }^2\Omega ^n\to \operatorname{Exp}%
(L^2(X;\sigma ))\otimes \bigg[
\bigoplus_{m=1}^nL_\sigma ^2\Psi ^n(X^m)\bigg]
\]
given by
\[
{\cal I}^n:=(U\otimes {\bf 1})I^n,
\]
where $U$ is the Wiener--It\^{o}--Segal isomorphism between the Poisson
space $L_{\pi _\sigma }^2(\Gamma _X)$ and the symmetric Fock space $%
\operatorname{Exp}(L^2(X;\sigma))$ over $L^2(X;\sigma )${\rm , } see e{\rm .}%
g{\rm .}\ {\rm \cite{AKR1}.}
\end{corollary}


\section{Dirichlet operators on differential forms over configuration spaces%
\label{dodf}}

In this section, we introduce Dirichlet operators associated with the
Poisson measure on $\Gamma _X$ which act in the spaces of square-integrable
forms. These operators generalize the notions of Bochner and de Rham--Witten
Laplacians on finite-dimensional manifolds. But first, we recall some known
facts and definitions concerning the usual Dirichlet operator of the Poisson
measure and Laplace operators on differential forms over finite-dimensional
manifolds.

\subsection{The intrinsic Dirichlet operator on functions}
\label{subsec3.1}

In this subsection, we recall some theorems from \cite{AKR1} which concern
the intrinsic Dirichlet operator in the space $L^2_{\pi_\sigma} (\Gamma_X)$,
to be used later.

Let us recall that the logarithmic derivative of the measure $\sigma $ is
given by the vector field
\begin{equation}
X\ni x\mapsto \beta _\sigma (x):=\frac{\nabla ^X\rho (x)}{\rho (x)}\in T_xX
\nonumber
\end{equation}
(where as usually $\beta _\sigma :=0$ on $\{\rho =0\}$). We wish
now to define a logarithmic derivative of the Poisson measure, and
for this we need a generalization of the notion of vector field.

For each $\gamma \in \Gamma _X$, consider the triple
\begin{equation}
T_{\gamma ,\,\infty }\Gamma _X\supset T_\gamma \Gamma _X\supset T_{\gamma
,0}\Gamma _X.  \nonumber
\end{equation}
Here, $T_{\gamma ,0}\Gamma _X$ consists of all finite sequences from $%
T_\gamma \Gamma _X$, and $T_{\gamma ,\,\infty }\Gamma _X%
\mbox{$:=$}%
\left( T_{\gamma ,0}\Gamma _X\right) ^{\prime }$ is the dual space, which
consists of all sequences $V(\gamma )=(V(\gamma )_x)_{x\in \gamma }$, where $%
V(\gamma )_x\in T_xX$. The pairing between any $V(\gamma )\in T_{\gamma
,\,\infty }\Gamma _X$ and $v(\gamma )\in T_{\gamma ,0}\Gamma _X$ with
respect to the zero space $T_\gamma \Gamma _x$ is given by
\begin{equation}
\langle V(\gamma ),v(\gamma )\rangle _\gamma =\sum_{x\in \gamma }\langle
V(\gamma )_x,v(\gamma )_x\rangle _x  \nonumber
\end{equation}
(the series is, in fact, finite). From now on, under a vector field over $%
\Gamma _X$ we will understand mappings of the form $\Gamma _X\ni \gamma
\mapsto V(\gamma )\in T_{\gamma ,\infty }\Gamma _X$.

The logarithmic derivative of the Poisson measure $\pi _\sigma $ is defined
as the vector field
\begin{equation}
\Gamma _X\ni \gamma \mapsto B_{\pi _\sigma }(\gamma )=(\beta _\sigma
(x))_{x\in \gamma }\in T_{\gamma ,\infty }\Gamma _X
\end{equation}
(i.e., the logarithmic derivative of the Poisson measure is the lifting of
the logarithmic derivative of the underlying measure).

The following theorem is a version of Theorem~3.1 in \cite{AKR1} (for more
general classes of functions and vector fields).

\begin{theorem}[Integration by parts formula on the Poisson space]
\label{th-ibp}$\text{{}}$\newline
For arbitrary $F^{(1)},F^{(2)}\in {\cal FC}$ and a smooth cylinder vector
field $V\in {\cal FV}$ $(:={\cal F}\Omega ^1)${\rm , } we have
\begin{gather*}
\int_{\Gamma _X}\nabla _V^\Gamma F^{(1)}(\gamma )F^{(2)}(\gamma )\,\pi
_\sigma (d\gamma )=-\int_{\Gamma _X}F^{(1)}(\gamma )\nabla _V^\Gamma
F^{(2)}(\gamma )\,\pi _\sigma (d\gamma ) \\
\text{{}}-\int_{\Gamma _X}F^{(1)}(\gamma )F^{(2)}(\gamma )\big[ \left\langle
B_{\pi _\sigma }(\gamma ),V(\gamma )\right\rangle _\gamma +\operatorname{div}%
^\Gamma V(\gamma )\big] \,\pi _\sigma (d\gamma ),
\end{gather*}
where the divergence $\operatorname{div}^\Gamma V(\gamma )$ of the vector
field $V$ is given by
\begin{gather*}
\operatorname{div}V(\gamma )=\sum_{x\in \gamma }\operatorname{div}%
_x^XV(\gamma )=\langle \operatorname{div}_{\bullet }^XV(\gamma ),\gamma
\rangle , \\
\operatorname{div}_x^XV(\gamma )%
\mbox{$:=$}%
\operatorname{div}^XV_x(\gamma ,x),\qquad x\in \gamma ,
\end{gather*}
$\operatorname{div}^X$ denoting the divergence on $X$ with respect to the
volume measure $m.$
\end{theorem}

\noindent
{\it Proof}. The theorem follows from formula (\ref{3.1}) and the usual
integration by parts formula on the space $L^2(\Lambda ^n,\sigma ^{\otimes
n})$ (see also the proof of Theorem~\ref{th4.1} below).\quad $\blacksquare $
\vspace{2mm}

Following \cite{AKR1}, we consider the intrinsic pre-Dirichlet
form on the Poisson space
\begin{equation}
{\cal E}_{\pi _\sigma }(F^{(1)},F^{(2)})=\int_{\Gamma _X}\langle \nabla
^\Gamma F^{(1)}(\gamma ),\nabla ^\Gamma F^{(2)}(\gamma )\rangle _\gamma
\,\pi _\sigma (d\gamma )  \label{3.2}
\end{equation}
with domain $D({\cal E}_{\pi _\sigma }):={\cal FC}$. By using the fact that
the measure $\pi _\sigma $ has all moments finite and noting that there
exists a function $\varphi\in{\cal D}$, $\varphi\ge0$, such that
\[
|\langle \nabla ^\Gamma F^{(1)}(\gamma ),\nabla ^\Gamma F^{(2)}(\gamma
)\rangle _\gamma|\le \langle \varphi,\gamma\rangle,
\]
one concludes that the expression (\ref{3.2}) is well-defined.

Let $H_\sigma $ denote the Dirichlet operator in the space $L^2(X;\sigma )$
associated to the pre-Dirichlet form
\begin{equation}
{\cal E}_\sigma (\varphi ,\psi )=\int_X\langle \nabla ^X\varphi (x),\nabla
^X\psi (x)\rangle _x\,\sigma (dx),\qquad \varphi ,\psi \in {\cal D}.
\nonumber
\end{equation}
This operator acts as follows:
\begin{equation}
H_\sigma \varphi (x)=-\Delta ^X\varphi (x)-\langle \beta _\sigma (x),\nabla
^X\varphi (x)\rangle _x,\qquad \varphi\in{\cal D},  \nonumber
\end{equation}
where $\Delta ^X:=\operatorname{div}^X\nabla ^X$ is the Laplace--Beltrami
operator on $X$.

Then, by using Theorem~\ref{th-ibp}, one gets
\begin{equation}
{\cal E}_{\pi _\sigma }(F^{(1)},F^{(2)})=\int_{\Gamma _X}H_{\pi _\sigma
}F^{(1)}(\gamma )F^{(2)}(\gamma )\,\pi _\sigma (d\gamma ),\qquad
F^{(1)},F^{(2)}\in {\cal FC}.  \label{sukk}
\end{equation}
Here, the intrinsic Dirichlet operator $H_{\pi _\sigma }$ is given by
\begin{align}
H_{\pi _\sigma }F(\gamma )
\mbox{$:=$}%
\sum_{x\in \gamma }H_{\sigma ,x}F(\gamma )\equiv \langle H_{\sigma ,\bullet
}F(\gamma ),\gamma \rangle ,  \nonumber \\
H_{\sigma ,x}F(\gamma )
\mbox{$:=$}%
H_\sigma F_x(\gamma ,x),\qquad x\in \gamma ,  \label{dir-op1}
\end{align}
so that the operator $H_{\pi _\sigma }$ is the lifting to $%
L^2_{\pi_\sigma}(\Gamma _X )$ of the operator $H_\sigma $ in $L^2(X;\sigma )$%
.

Upon (\ref{sukk}), the pre-Dirichlet form ${\cal E}_{\pi _\sigma }$ is
closable, and we preserve the notation for the closure of this form.

\begin{theorem}
\label{th3.2}{\rm \cite{AKR1}} Suppose that $(H_\sigma ,{\cal D})$
is essentially self-adjoint on $L^2(X;\sigma )${\rm . } Then{\rm ,
} the operator $H_{\pi _\sigma }$ is essentially self-adjoint on
${\cal FC}.$
\end{theorem}

\begin{remark}
\label{proof}{\rm This theorem was proved in \cite{AKR1}, Theorem~5.3. (We
have already mentioned in Remark~\ref{rem2.1} that the inclusion ${\cal FC}%
_{\mathrm b}^\infty ({\cal D},\Gamma _X)\subset {\cal FC}$ holds.)
We would like to stress that this result is based upon the theorem
which says that the image of the operator $H_{\pi _\sigma }$ under
the isomorphism $U$ between the Poisson space and the Fock space
$\operatorname{Exp}\left( L^2(X;\sigma
)\right) $ is the differential second quantization $d\operatorname{Exp}%
H_\sigma $ of the operator $H_\sigma $. }
\end{remark}

\begin{remark}\label{rem_Local}
{\rm The condition of Theorem~\ref{th3.2} is satisfied if e.g.\ }
\begin{equation}
\Vert \beta _\sigma \Vert _{TX}\in L_{\mathrm loc}^p(X;\sigma )
\label{est2}
\end{equation}
{\rm \ for some $p>\dim \,X,$ see \cite{AKR1}.}
\end{remark}

In what follows, we will suppose for simplicity that
\begin{equation}\label{cgg}\text{the function $\rho $ is infinitely differentiable on
$X$ and $\rho (x)>0$ for all $x\in X$.}\end{equation} Evidently,
estimate (\ref{est2}) is implied by \eqref{cgg}.

Finally, we mention the important fact  \cite{AKR1} that the
diffusion process which is
properly associated with the Dirichlet form $({\cal E}_{\pi _\sigma },D(%
{\cal E}_{\pi _\sigma }))$ is the usual independent infinite
particle process (or distorted Brownian motion on $\Gamma _X$),
introduced by Doob \cite{Doob} .

\subsection{Laplacians on differential forms over finite-dimensional
manifolds}

We recall now some facts on the Bochner and de Rham--Witten Laplacians on
differential forms over a finite-dimensional manifold.

Let $M$ be a Riemannian manifold equipped with the measure $\mu (dx)=e^{\phi
(x)}dx,$ $dx$ being the volume measure and $\phi $ a $C^2$-function on $M$.
We consider a Hilbert bundle
\begin{equation}
{\cal H}_x\mapsto x\in M  \nonumber
\end{equation}
over $M$ equipped with a smooth connection, and denote by $\nabla $ the
corresponding covariant derivative in the spaces of sections of this bundle.
Let $L^2(M\rightarrow {\cal H};\mu )$ be the space of $\mu $-square
integrable sections. The operator
\begin{equation}H_\mu^{\mathrm B}:=
\nabla _\mu ^{*}\nabla  \nonumber
\end{equation}
in $L^2(M\rightarrow {\cal H};\mu )$, where $\nabla _\mu ^{*}$ is
the adjoint of $\nabla $, will be called the Bochner Laplacian
associated with the measure $\mu $.  One can easily write the
corresponding differential expression on the space of twice
differentiable sections. In the case where $\phi \equiv 0$ and ${\cal H}%
_x=\wedge ^n(T_xM)$, we obtain the classical Bochner Laplacian on
differential forms (see e.g.\ \cite{CFKSi}).

Now, let $d$ be the exterior differential in spaces of differential forms
over $M.$ The operator
\begin{equation}
H_\mu ^{\mathrm R}%
\mbox{$:=$}%
d_\mu ^{*}d+dd_\mu ^{*}  \nonumber
\end{equation}
acting in the space of $\mu $-square integrable forms, where
$d_\mu ^{*}$ is the adjoint of $d$, will be called the de Rham
Laplacian associated with the measure $\mu $ (or the Witten
Laplacian associated with $\phi $, see e.g.\  \cite {CFKSi}).

We will use sometimes more extended notations $H_{\mu ,n}^{\mathrm
B}(M),\;H_{\mu ,n}^{\mathrm R}(M)$ for the Bochner and de
Rham--Witten Laplacians on the space of $\mu$-square integrable
$n$-forms over $M$.

The relation of the Bochner and de Rham--Witten Laplacians on differential
forms is given by the Weitzenb\"{o}ck formula (cf.\ \cite{CFKSi}, \cite{E3}%
), which will be recalled now.

Fix $x\in M$ and let $(e_j)_{j=1}^{\operatorname{dim}M}$ be an orthonormal
basis in $T_xM$. Denote by
\begin{align}
a_j&:\wedge^{n+1}(T_xM)\to\wedge^n(T_xM),  \nonumber \\
a_j^*&:\wedge^{n}(T_xM)\to\wedge^{n+1}(T_xM)  \label{odyn}
\end{align}
the annihilation and creation operators, respectively, defined by
\begin{alignat}{3}
a_ju^{(n+1)}&=\sqrt{n+1}\,\langle e_j,u^{(n+1)}\rangle_x, &
&\qquad & &u^{(n+1)}\in\wedge^{n+1}(T_xM),  \nonumber \\
a_j^*u^{(n)}&=\sqrt{n+1}\,e_j\wedge u^{(n)}, & &\qquad & &u^{(n)}\in%
\wedge^n(T_xM)  \label{dvi}
\end{alignat}
(the pairing in the expression $\langle e_j,u^{(n+1)}\rangle_x$ is
carried out in the first ``variable,'' so that $a^*_j$ becomes
adjoint of $a_j$).

Let us introduce the operator $R_n(x)$ in $\wedge^n(T_xM)$ by
\[
R_n(x):=\sum_{i,j,k,l=1}^{\operatorname{dim}M}
R_{ijkl}(x)\,a_i^*a_ja_k^*a_l,
\]
where $R_{ijkl}$ is the curvature tensor on $M$. It can be shown that the
definition of this operator is independent of the specific choice of basis.

Next, let $(\nabla^M\beta_\mu(x))^{\wedge n}$ be the operator in $%
\wedge^n(T_xM)$ given by
\begin{multline}  \label{6.1}
(\nabla^M\beta_\mu(x))^{\wedge n}:=\nabla^M\beta_\mu(x) \otimes{\bf 1}%
\dots\otimes {\bf 1}+ {\bf 1}\otimes \nabla^M\beta_\mu(x)\otimes {\bf 1}%
\otimes\dots\otimes{\bf 1} \\
\text{}+\dots+{\bf 1}\otimes\dots\otimes{\bf 1}\otimes \nabla^M\beta_\mu(x),
\end{multline}
$\nabla^M\beta_\mu(x)$ being understood as an operator in $T_xM$.

Then, the Weitzenb\"{o}ck formula writes as follows:
\begin{equation}
H_\mu ^{\mathrm R}\omega_n(x)=H_\mu ^{\mathrm B}\omega_n(x)+R_\mu
(x)\omega_n(x), \label{wei1}
\end{equation}
where $\omega_n$ is an n-form on $X$, and $R_\mu
(x)\omega_n(x)=R_{\mu ,n}(x)\omega_n(x)$,
\begin{equation}
R_{\mu ,n}(x):=R_n(x)-(\nabla ^M\beta _\mu (x))^{\wedge n}.  \label{wei2}
\end{equation}

\begin{remark}\rom{The classical Weitzenb\"{o}ck formula
is related, in fact, to the case where $\phi \equiv 0$, see e.g.\
\cite{CFKSi}, \cite{E3}. Formula (\ref{wei1}) can be obtained by a
direct calculation using similar arguments, cf.\ \cite{ADK1}.}
\end{remark}

\subsection{Bochner Laplacian on forms over the configuration space}

Let us consider the pre-Dirichlet form
\begin{equation}  \label{4.2}
{\cal E}^{\mathrm B}_{\pi_\sigma}(W^{(1)},W^{(2)})=\int_{\Gamma_X}
\langle\nabla^\Gamma W ^{(1)}(\gamma),\nabla^\Gamma W
^{(2)}(\gamma)\rangle_{T_\gamma\Gamma_X
\otimes\wedge^n(T_\gamma\Gamma_X)}\,\pi_\sigma(d\gamma),
\end{equation}
where $W^{(1)},W^{(2)}\in {\cal F}\Omega^n$. As easily seen, there
exists $\varphi\in{\cal D}$, $\varphi\ge0$, such that
\[
|\langle\nabla^\Gamma W ^{(1)}(\gamma),\nabla^\Gamma
W^{(2)}(\gamma)\rangle_{T_\gamma\Gamma_X
\otimes\wedge^n(T_\gamma\Gamma_X)}|\le \langle
\varphi^{\otimes(n+1)},\gamma^{\otimes (n+1)}\rangle,
\]
so that the function under the sign of integral in \eqref{4.2} is integrable
with respect to $\pi_\sigma$.

\begin{theorem}
\label{th4.1} For any $W^{(1)},W^{(2)}\in {\cal F}\Omega ^n$, we have
\[
{\cal E}_{\pi _\sigma }^{\mathrm B}(W^{(1)},W^{(2)})=\int_{\Gamma
_X}\langle H_{\pi _\sigma }^{\mathrm B}W^{(1)}(\gamma
),W^{(2)}(\gamma )\rangle _{\wedge ^n(T_\gamma \Gamma _X)}\,\pi
_\sigma (d\gamma ),
\]
where $H_{\pi _\sigma }^{\mathrm B}$ is the operator in the space
$L_{\pi _\sigma }^2\Omega ^n$ with domain ${\cal F}\Omega ^n$
given by
\begin{equation}
H_{\pi _\sigma }^{\mathrm B}W(\gamma )=-\Delta ^\Gamma W(\gamma
)-\langle \nabla ^\Gamma W(\gamma ),B_{\pi _\sigma }(\gamma
)\rangle _\gamma ,\qquad W\in {\cal F}\Omega ^n.  \label{boch1}
\end{equation}
Here{\rm , }
\begin{equation}
\Delta ^\Gamma W(\gamma )%
\mbox{$:=$}%
\sum_{x\in \gamma }\Delta _x^XW(\gamma )\equiv \left\langle \Delta _{\bullet
}^\Gamma W(\gamma ),\gamma \right\rangle ,  \label{boch2}
\end{equation}
where $\Delta _x^X$ is the Bochner Laplacian of the bundle $\wedge
^n(T_{\gamma _y}\Gamma _X)\mapsto y\in {\cal O}_{\gamma ,x}$ with the volume
measure{\rm . }
\end{theorem}

\noindent{\it Proof}. First, we note that, for $W\in{\cal F}\Omega^n$,
\[
\Delta^X_xW(\gamma)_{[x_1,\dots,x_n]_d}:=\Delta^XW_x(\gamma,x)_
{[x_1,\dots,x_n]_d},\qquad x\in\gamma,\,[x_1,\dots,x_n]_d\subset\gamma.
\]

Fix now $W^{(1)},W^{(2)}\in{\cal F}\Omega^n$ and let $\Lambda_1, \Lambda_2$
be compact subsets of $X$ as in Definition~\ref{def2.2} corresponding to $%
W^{(1)}$, $W^{(2)}$, respectively. Let $\Lambda$ be an open set in $X$ with
compact closure such that both $\Lambda_1$ and $\Lambda_2$ are subsets of $%
\Lambda$. Then, by using \eqref{3.1},
\begin{gather*}
\int_{\Gamma_X}\langle\nabla^\Gamma W^{(1)}(\gamma),\nabla^\Gamma
W^{(2)}(\gamma)\rangle
_{T_\gamma\Gamma_X\otimes\wedge^n(T_\gamma\Gamma_X)}\,\pi_\sigma(d\gamma)=
\\ = e^{-\sigma(\Lambda)}\sum_{k=0}^\infty\frac1{k!}
\int_{\Lambda^k}\sum_{i=1}
^k\langle\nabla^X_{x_i}W^{(1)}(\{x_1,\dots,x_k\}), \\
\nabla^X_{x_i}W^{(2)}(\{x_1,\dots,x_k\})
\rangle_{T_{x_i}X\otimes\wedge^n(T_{x_1}X\oplus\dots\oplus
T_{x_k}X )}\,\sigma(dx_1)\dotsm\sigma(dx_k) \\
=e^{-\sigma(\Lambda)}\sum_{k=0}^\infty\frac1{k!}\sum_{i=1}^k\int_{\Lambda^k}
\sum_{[y_1,\dots,y_n]_d\subset \{x_1,\dots,x_k\}}
\langle\nabla^X_{x_i}W^{(1)}(\{x_1,\dots,x_k\})
_{[y_1,\dots,y_n]_d}, \\ \nabla^X_{x_i}W^{(2)}(\{x_1,\dots,x_k\})
_{[y_1,\dots,y_n]_d}\rangle_{T_{x_i}X\otimes(T_{y_1}X\wedge\dots\wedge
T_{y_n}X)}\, \sigma(dx_1)\dotsm\sigma(dx_k) \\
=e^{-\sigma(\Lambda)}\sum_{k=0}^\infty\frac1{k!}\sum_{i=1}^k\int_{\Lambda^k}
\sum_{[y_1,\dots,y_n]_d\subset \{x_1,\dots,x_k\}} \langle
\Delta^X_{x_i}W^{(1)}(\{x_1,\dots, x_k\})_{[y_1,\dots,y_n]_d} \\
\text{}+\langle\nabla^X_{x_i}W^{(1)}(\{x_1,\dots,
x_k\})_{[y_1,\dots,y_n]_d},\beta_\sigma(x_i)\rangle_{x_i}, \\
W^{(2)}(\{x_1,\dots,x_k\}) _{[y_1,\dots,y_n]_d}\rangle_{
T_{y_1}X\wedge\dots\wedge T_{y_n}X}\,
\sigma(dx_1)\dotsm\sigma(dx_k) \\ =\int_{\Gamma_X}\langle
H^{\mathrm
B}_{\pi_\sigma}W^{(1)}(\gamma),W^{(2)}(\gamma)\rangle_{\wedge^n
(T_\gamma\Gamma_X)}\,\pi_\sigma(d\gamma).\quad\blacksquare
\end{gather*}

\begin{remark}
\label{rem4.1}{\rm We can rewrite the action of the operator
$H_{\pi _\sigma }^{\mathrm B}$ in the two following forms: }

\begin{enumerate}
\item[{\rm 1)}]  {\rm We have from (\ref{boch1}) and (\ref{boch2}) that
\begin{equation}
H_{\pi _\sigma }^{\mathrm B}W(\gamma )=\sum_{x\in \gamma
}H_{\sigma ,x}^{\mathrm B}W(\gamma )\equiv \left\langle H_{\sigma
,\bullet }^{\mathrm B}W(\gamma ),\gamma \right\rangle ,\qquad
W(\gamma )\in {\cal F}\Omega ^n, \label{blo1}
\end{equation}
where
\begin{equation}
H_{\sigma ,x}^{\mathrm B}W(\gamma ):=-\Delta _x^XW(\gamma
)-\left\langle \nabla _x^XW(\gamma ),\beta _\sigma
(x)\right\rangle _x. \label{blo2}
\end{equation}
Thus, the operator $H_{\pi_\sigma}^{\mathrm B}$ is the lifting of
the Bochner Laplacian on $X$ with the measure $\sigma .$ }

\item[{\rm 2)}]  {\rm As easily seen, the operator $H_{\pi _\sigma }^{\mathrm B}$ preserves
the space ${\cal F}\Omega ^n$, and we can always take $\Lambda
(H_{\pi _\sigma }^{\mathrm B}W)=\Lambda (W)$. Then, for any open
bounded $\Lambda \supset \Lambda (W)$ (cf.\
Remark~\ref{form-fin}), we have
\begin{equation}
(H_{\pi _\sigma }^{\mathrm B}W)_{\Lambda ,\gamma }=H_{\sigma
^{\otimes |\Lambda \cap \gamma |}}^{\mathrm B}(X^{|\Lambda \cap
\gamma |})W_{\Lambda ,\gamma }, \label{cyl-boch}
\end{equation}
where }$H_{\sigma ^{\otimes |\Lambda \cap \gamma |}}^{\mathrm
B}(X^{|\Lambda \cap \gamma |})$ {\rm  is the Bochner Laplacian of
the manifold $X^{|\Lambda \cap
\gamma |}$
 with the product
measure $\sigma ^{\otimes |\Lambda \cap \gamma |}$ (cf.\
(\ref{cyl-form})). The equality \eqref{cyl-boch} holds on ${\cal
O} _{\gamma,x_1}\times\dots\times {\cal
O}_{\gamma,x_{|\Lambda\cap\gamma|}}$, where
$\{x_1,\dots,x_{|\Lambda\cap\gamma|}\}=\Lambda\cap\gamma$. Notice
that, since the operator $H_{\sigma ^{\otimes |\Lambda \cap \gamma
|}}^{\mathrm B}(X^{|\Lambda \cap \gamma |})$ acts locally on
(smooth) forms on $X^{|\Lambda\cap\gamma|}$, the expression  on
the right hand side of \eqref{cyl-boch} is well defined as a form
on ${\cal O} _{\gamma,x_1}\times\dots\times {\cal
O}_{\gamma,x_{|\Lambda\cap\gamma|}}$. }
\end{enumerate}
\end{remark}

It follows from Theorem~\ref{th4.1} that the pre-Dirichlet form ${\cal E}%
_{\pi _\sigma }^{\mathrm B}$ is closable in the space $L_{\pi
_\sigma }^2\Omega ^n$. The generator of its closure (being
actually the Friedrichs extension of the operator $H_{\pi _\sigma
}^{\mathrm B}$, for which we will use the same notation) will be
called the Bochner Laplacian on $n$-forms over $\Gamma _X$
corresponding to the Poisson measure $\pi _\sigma $.

For linear operators $A$ and $B$ acting in Hilbert spaces ${\cal H}$ and $%
{\cal K}$, respectively, we introduce the operator $A\boxplus B$ in ${\cal H}%
\otimes {\cal K}$ by
\begin{equation}
A\boxplus B%
\mbox{$:=$}%
A\otimes {\bf 1}+{\bf 1}\otimes B,\qquad \operatorname{Dom}(A\boxplus B):=%
\operatorname{Dom}(A)\otimes _{\mathrm a}\operatorname{Dom}(B),
\nonumber
\end{equation}
where $\otimes _{\mathrm a}$ stands for the algebraic tensor
product. Next, for
operators $A_1,\dots ,A_n$ acting in Hilbert spaces ${\cal H}_1,\dots ,{\cal %
H}_n$, respectively, let $\bigoplus_{i=1}^nA_i$ denote the operator in $%
\bigoplus_{i=1}^nH_i$ given by
\[
\bigg(\bigoplus_{i=1}^nA_i\bigg)(f_1,\dots ,f_n)=(A_1f_1,\dots
,A_nf_n),\qquad f_i\in \operatorname{Dom}(A_i).
\]

\begin{theorem}
\label{thonsa1} {\rm 1)} On ${\cal D}\Omega ^n$ we have
\begin{equation}
H_{\pi _\sigma }^{\mathrm B}=(I^n)^{-1}\bigg[H_{\pi _\sigma }\boxplus \bigg( %
\bigoplus_{m=1}^nH_{\sigma ,\,(n,m)}^{\mathrm B}\bigg)\bigg]I^n,
\label{dec-gen0}
\end{equation}
where $H_{\sigma ,\,(n,m)}^{\mathrm B}$ denotes the restriction of
the Bochner Laplacian $H_{\sigma ^{\otimes m},n}^{\mathrm B}(X^m)$
acting in the space $L^2(X^m\to \wedge ^n(TX^m);\sigma ^{\otimes
m})$ to the subspace $L_\sigma ^2\Psi ^n(X^m)${\rm . }

{\rm 2)} Suppose that, for each $m=1,\dots ,n$\rom, the Bochner
Laplacian $H_{\sigma ,m}^{\mathrm B}(X)$ is essentially
self-adjoint on the set of smooth forms with compact support.
Then{\rm , } ${\cal D}\Omega ^n$ is a domain of
essential self-adjointness of $H_{\pi _\sigma }^B${\rm , } and the equality~%
\eqref{dec-gen0} holds for the closed operators $H_{\pi _\sigma }^{\mathrm B}$ and $%
H_{\pi _\sigma }\boxplus \big( \bigoplus_{m=1}^nH_{\sigma
,\,(n,m)}^{\mathrm B}\big)$ {\rm (}where the latter
operator is closed from its domain of essential self-adjointness $I^n({\cal D}%
\Omega ^n)${\rm )}{\rm . }
\end{theorem}

\begin{remark}\label{brand_new_remark1}\rom{ The essential
self-adjointness of the Bochner Laplacian $H^{\mathrm B}_{\sigma}$
on the set of smooth forms with compact support is well-known in
the case where $\sigma$ is the volume measure, see e.g.\
\cite{E3}. More generally, it is sufficient to assume that
$\beta_\sigma$, together with its derivatives up to order 2, is
bounded.}\end{remark}

\noindent {\it Proof of Theorem}~\ref{thonsa1}. 1) Let $W\in{\cal D}\Omega^n$ be given by the formula %
\eqref{kikimora}. Then, using \eqref{blo1}, \eqref{blo2}, and \eqref{dir-op1}%
, we get
\begin{gather}
(H^{B}_{\pi_\sigma}W)_k(\gamma)=0\qquad \text{for }k\ne m,  \nonumber \\
(H^{B}_{\pi_\sigma}W)_m(\gamma)(\bar{x})= \bigg(\sum_{x\in%
\gamma}H^{B}_{\sigma,x}W\bigg)_m(\gamma)(\bar {x}) \nonumber
\\
=\bigg(\sum_{x\in\gamma\setminus\{\bar{x}\}}H^{B}_{\sigma,x}W\bigg)%
_m(\gamma)(\bar {x}) +\bigg(\sum_{x\in\{\bar{x}%
\}}H^{B}_{\sigma,x}W\bigg)_m(\gamma)(\bar{x})  \nonumber
\\
=(m!)^{1/2}\,\left[\bigg(\sum_{x\in\gamma\setminus\{\bar{x}\}}H_{\sigma,x}F\bigg)%
(\gamma\setminus \{\bar{x}\})\omega(\bar{x})+F(\gamma\setminus
\{\bar{x}\})\bigg(\sum_{x\in\{\bar{x}\}}H^{B}_{\sigma,x}\omega%
\bigg)(\bar{x})\right]  \nonumber \\
=(m!)^{1/2}\,\left[(H_{\pi_\sigma}F)(\gamma\setminus \{\bar{x}\})\omega(\bar{x}%
)+F(\gamma\setminus \{\bar{x}\})(H^{B}_{\sigma,\,(n,m)}\omega)(%
\bar{x})\right ].  \label{chyzh}
\end{gather}
(Notice that the Bochner Laplacian in the space $L^2(X^m\to\wedge^n(TX^m);%
\sigma^{\otimes m})$ leaves the set $\Psi_0^n(X^m)$ invariant.) Therefore,
\begin{equation}  \label{pyzh}
(I_k^n H^{B}_{\pi_\sigma}W)(\gamma,\bar x)=\begin{cases}0,&
\text{for }k\ne m,\\ (H_{\pi_\sigma}F)
(\gamma)\omega(\bar{x})+F(\gamma)(H^{\mathrm B
}_{\sigma,\,(n,m)}\omega)(\bar{x}),&\text{for }k=m.\end{cases}
\end{equation}
Hence, by virtue of \eqref{is1}, we get
\[
\bigg(\bigg[H_{\pi_\sigma}\boxplus\bigg(\bigoplus_{i=1}^n H^{B}_{\sigma,(n,i)}%
\bigg)\bigg]I^nW\bigg) _k(\gamma,\bar x)=(I_k^n
H^{B}_{\pi_\sigma}W)(\gamma,\bar x),\qquad k=1,\dots,n
\]
which proves \eqref{dec-gen0}.

2) Let $\Omega _0^n(X^m)$ denote the set of all smooth forms
$\omega \colon X^m\to \wedge ^n(X^m)$ with compact support. It is
not hard to see that  the essential self-adjointness of $H_{\sigma
,m}^{\mathrm B}(X) $ for each $m=1,\dots ,n$ implies that
\begin{equation}\label{statement}
\text{ the Bochner Laplacian $H^{\mathrm B }:=H_{\sigma ^{\otimes
m},n}^{\mathrm B}(X^m)$ is essentially self-adjoint on $\Omega
_0^n(X^m)$.}\end{equation}

 Indeed, by using the decomposition
\eqref{n-forms}, we have \begin{gather*}
L^2(X^m\to\wedge^n(TX^m);\sigma^{\otimes m})=L^2( X^m\ni
(x_1,\dots,x_m)\to
\wedge^n(T_{(x_1,\dots,x_m)}X^m);\sigma^{\otimes m})\\=
\bigoplus_{\begin{gathered}{\scriptstyle{
0\le k_1,\dots,k_m\le d}} \\ \scriptstyle k_1+\dots+k_m=n \end{gathered}%
}L^2( X^m\ni(x_1,\dots,x_m)\to(T_{x_1}X)^{\wedge
k_1}\wedge\dots\wedge(T_{x_m}X)^{\wedge k_m};\sigma^{\otimes
m}),\end{gather*} and it is enough to show that the Bochner
Laplacian $H ^{\mathrm B}$ is essentially self-adjoint  in each
space
\begin{equation}\label{eins}L^2(X^m\ni(x_1,\dots,x_m)\to(T_{x_1}X)^{\wedge
k_1}\wedge\dots\wedge(T_{x_m}X)^{\wedge k_m};\sigma^{\otimes
m})\end{equation} on the set of smooth forms.

On the other hand, by using the essential self-adjointness of each
operator $H^{\mathrm B}_{\sigma,m}(X)$ on $\Omega_0^m(X)$ and that
of the operator $H_\sigma$ in the space $L^2(X;\sigma)$ on the set
$\cal D$ (Remark~\ref{rem_Local}), we  conclude from the theory of
operators admitting separation of variables \cite[Ch.~6]{B} that
the operator \begin{equation}\label{drei} H^{\mathrm
B}_{\sigma,k_1}(X)\sqplus\cdots\sqplus H^{\mathrm
B}_{\sigma,k_m}(X),\qquad H^{\mathrm
B}_{\sigma,0}(X):=H_\sigma,\end{equation} is essentially
self-adjoint in the space
\begin{equation}\label{zwei}\begin{gathered}
L^2(X\to\wedge^{k_1}(TX);\sigma)\otimes\dots\otimes
L^2(X\to\wedge^{k_m}(TX);\sigma)=\\=L^2(X^m\ni(x_1,\dots,x_m)\to
(T_{x_1}X)^{\wedge k_1}\otimes\dots\otimes(T_{x_m}X)^{\wedge
k_m};\sigma^{\otimes m})\end{gathered}\end{equation} on the
algebraic product of the domains of the operators $H^{\mathrm
B}_{\sigma,k_i}(X)$.

Next, we note that, for each $(x_1,\dots,x_m)\in\tilde X^m$, there
exists an intrinsic unitary isomorphism
$$\operatorname{Iso}_{k_1,\dots,k_m}\colon (T_{x_1}X)^{\wedge
k_1}\otimes\dots\otimes(T_{x_m}X)^{\wedge
k_m}\to(T_{x_1}X)^{\wedge k_1}\wedge\dots\wedge(T_{x_m}X)^{\wedge
k_m}$$ that is given by the formula
\begin{multline*}\operatorname{Iso}_{k_1,\dots,k_m}(u_1^{(1)}\wedge\dots\wedge
u_{k_1}^{(1)})\otimes \dots\otimes (u_1^{(m)}\wedge\dots\wedge
u_{k_m}^{(m)}):=\\ =\sqrt{\dfrac{(k_1+\dots+k_m)!}{k_1!\dotsm
k_m!}}\,u_1^{(1)}\wedge\dots\wedge u_{k_1}^{(1)}\wedge \dots\wedge
u_1^{(m)}\wedge\dots\wedge u_{k_m}^{(m)},\qquad u_j^{(i)}\in
T_{x_i}X,\end{multline*}
 and then it is extended by
linearity. As easily seen, this definition is independent of the
representation of a vector from $(T_{x_1}X)^{\wedge
k_1}\otimes\dots\otimes(T_{x_m}X)^{\wedge k_m}$.  Hence, for any
$(k_1,\dots,k_m)$, we can construct the unitary ${\cal
U}_{k_1,\dots,k_m}$ between the spaces \eqref{zwei} and
\eqref{eins} by setting $$({\cal
U}_{k_1,\dots,k_m}F)(x_1,\dots,x_m):=\operatorname{Iso}_{k_1,\dots,k_m}
(F(x_1,\dots,x_m)).$$ Under this unitary, the operator
\eqref{drei} goes over into the operator $H^{\mathrm B}$ in the
space \eqref{eins}, while the image of its domain consists of
linear combinations of the form ${\cal
U}_{k_1,\dots,k_m}(\omega^{(k_1)}\otimes\dots\otimes\omega^{(k_m)})$,
$\omega^{(k_i)}\in\Omega_0^{k_i}(X)$. From here, the assertion
\eqref{statement} follows.

Let $\widehat{L}^2(X^m\to \wedge ^n(TX^m);\sigma ^{\otimes m})$
denote the subspace of $L^2(X^m\to \wedge ^n(TX^m);\sigma
^{\otimes m})$ consisting of all symmetric forms, i.e., the forms
$\omega \in L^2(X^m\to \wedge ^n(TX^m);\sigma ^{\otimes m})$ for
which the equality \eqref{symmetric} holds for $\sigma^{\otimes m
}$-a.a.\ $(x_1,\dots,x_m)\in X^m$. Evidently, the orthogonal
projection $P_m^n$ onto this subspace is given by the formula
\begin{equation}
(P_m^n\omega )(x_1,\dots ,x_m)=\frac 1{m!}\sum_{\sigma \in
S_m}\omega (x_{\sigma (1)},\dots ,x_{\sigma (m)})  \label{kkk1}
\end{equation}
and
\begin{equation}
P_m^n\Omega _0^n(X^m)=\Omega _{0,\,\mathrm sym}^n(X^m),
\label{kkk2}
\end{equation}
where  $\Omega _{0,\,\mathrm sym}^n(X^m)$ denotes the set of
symmetric smooth forms $\omega \colon X^m\to \wedge ^n(X^m)$ with
compact support.

The assertion \eqref{statement} and the nonnegative definiteness
of $H^{\mathrm B}$ yield that the set $(H^{\mathrm B}+{\bf
1})\Omega _0^n(X^m)$ is dense in $L^2(X^m\to \wedge
^n(TX^m);\sigma ^{\otimes m})$, see e.g.\ \cite{RS}, Section~10.1.
Therefore, the set $P_m^n(H^{\mathrm B}+{\bf 1})\Omega _0^n(X^m)$ is dense in $%
\widehat{L}^2(X^m\to \wedge ^n(TX^m);\sigma ^{\otimes m})$. But upon %
\eqref{kkk1} and \eqref{kkk2},
\[
P_m^n(H^{\mathrm B}+{\bf 1})\Omega _0^n(X^m)=(H^{\mathrm B}P_m^n+P_m^n)\Omega _0^n(X^m)=
(H^{\mathrm B}+%
{\bf 1})\Omega _{0,\,{\mathrm sym}}^n(X^m),
\]
which implies that the Bochner Laplacian $H^{\mathrm B}$ in the space $\widehat{L}%
^2(X^m\to \wedge ^n(TX^m);\sigma ^{\otimes m})$ is essentially
self-adjoint on $\Omega _{0,\,\mathrm sym}^n(X^m)$.

Because $H^{\mathrm B}$ acts invariantly on the subspace $L_\sigma
^2\Psi ^n(X^m)$ and also on its orthogonal complement in $\widehat{L}%
^2(X^m\to \wedge ^n(TX^m);\sigma ^{\otimes m})$, we conclude that
$H_{\sigma ,\,(n,m)}^{\mathrm B}$ is essentially self-adjoint on
$\Psi _0^n(X^m)$. Consequently, the operator
$\bigoplus_{m=1}^nH_{\sigma ,\,(n,m)}^{\mathrm B}$ is essentially
self-adjoint on the direct sum of the sets $\Psi _0^n(X^m)$,
$m=1,\dots ,n$.

Finally, taking to notice that the operator $%
H_{\pi _\sigma }$ is essentially self-adjoint on ${\cal FC}$
(Theorem~\ref {th3.2}), we conclude again from the theory of
operators admitting separation of variables  that $I^n({\cal
D}\Omega ^n)$ is a domain of
essential self-adjointness of the operator $H_{\pi _\sigma }\boxplus \big( %
\bigoplus_{m=1}^nH_{\sigma ,\,(n,m)}^{\mathrm B}\big)$ in the
space $L_{\pi _\sigma }^2(\Gamma _X)\otimes \big[
\bigoplus_{m=1}^nL_\sigma ^2\Psi ^n(X^m)\big]$. Thus,
(\ref{dec-gen0}) yields the statement.\quad $\blacksquare
$\vspace{2mm}

We give also a Fock space representation of the operator $H_{\pi
_\sigma }^{\mathrm B} $. Corollary~\ref{fock} implies the
following

\begin{corollary}
\label{fock1} Let the conditions of Theorem~{\rm \ref{thonsa1}, 2)} be
satisfied{\rm . } Then{\rm , }
\begin{equation}
{\cal I}^nH_{\pi _\sigma }^{\mathrm B}({\cal
I}^n)^{-1}=d\operatorname{Exp}H_\sigma \boxplus \bigg(
\bigoplus_{m=1}^nH_{\sigma ,\,(n,m)}^{\mathrm B}\bigg),  \nonumber
\end{equation}
cf{\rm .}\ Remark~{\rm \ref{proof}.}
\end{corollary}


\subsection{De Rham Laplacian on forms over the configuration space}

We define linear operators
\begin{equation}  \label{dunia}
d^\Gamma\colon {\cal F}\Omega^n\to {\cal F}\Omega^{n+1},\qquad n\in{\Bbb N}%
_0,\quad {\cal F}\Omega^{0}:={\cal FC},
\end{equation}
by
\begin{equation}  \label{5.1}
(d^\Gamma W )(\gamma):=\sqrt{n+1}\,\operatorname{AS}_{n+1}(\nabla^\Gamma
W(\gamma)),
\end{equation}
where $\operatorname{AS}_{n+1}\colon (T_\gamma\Gamma_X)^{\otimes(n+1)}\to%
\wedge^{n+1}(T_\gamma\Gamma_X)$ is the antisymmetrization operator. It
follows from this definition that
\begin{equation}  \label{5.2}
(d^\Gamma W )(\gamma)=\sum_{x\in\gamma}(d^X_xW)(\gamma),
\end{equation}
where
\begin{align}
(d^X_xW)(\gamma):&=\sum_{[x_1,\dots,x_n]_d
\subset\gamma}d^X(W_x(\gamma,x)_{[x_1,\dots,x_n]_d})  \nonumber \\
&=\sum_{[x_1,\dots,x_n]_d\subset\gamma}\sqrt{n+1}\,\operatorname{AS}%
_{n+1}(\nabla^ XW_x(\gamma,x)_{[x_1,\dots, x_n]_d})  \label{5.3}
\end{align}
with $\operatorname{AS}_{n+1}\colon T_xX\otimes (T_{x_1}X\wedge\dots\wedge
T_{x_n}X)\to T_xX\wedge T_{x_1}X\wedge \dots\wedge T_{x_n}X$ being again the
antisymmetrization. Therefore, we have indeed the inclusion $%
d^\Gamma\omega\in{\cal F}\Omega^{n+1}$ for each $\omega\in{\cal F}\Omega^n$.

Suppose that, in local coordinates on the manifold $X$, the form $%
W_x(\gamma,\bullet)_{[x_1 ,\dots,x_n]_d}$ has the representation
\begin{equation}  \label{pils1}
{\cal O}_{\gamma,x}\ni y\mapsto W_x(\gamma,y)_
{[x_1,\dots,x_n]_d}=w(y)\,h_1\wedge\dots\wedge h_n,.
\end{equation}
where $w\colon{\cal O}_{\gamma,x}\to{\Bbb R}$ and $h_1\wedge\dots\wedge
h_n\in T_{x_1}\wedge\dots\wedge T_{x_n}$. Then,
\begin{equation}  \label{pils2}
\operatorname{AS}_{n+1}(\nabla^XW_x(\gamma,x)_{[x_1,\dots,x_n]_d})=\nabla^X
w(x) \wedge h_1\wedge\dots\wedge h_n,
\end{equation}
which, upon \eqref{5.3}, describes the action of $d_x^X$.

Let us consider now $d^\Gamma $ as an operator acting from the space $L_{\pi
_\sigma }^2\Omega ^n$ into $L_{\pi _\sigma }^2\Omega ^{n+1}$. Analogously to
the proof of Theorem~\ref{th4.1}, we get the following formula for the
adjoint operator $d_{\pi _\sigma }^{\Gamma *}$ restricted to ${\cal F}\Omega
^{n+1}$:
\begin{equation}
(d_{\pi _\sigma }^{\Gamma *}W)(\gamma )=\sum_{x\in \gamma }(d_{\sigma
,x}^{X*}W)(\gamma ),  \label{5.4}
\end{equation}
where
\begin{equation}
(d_{\sigma ,x}^{X*})W(\gamma )=\sum_{[x_1,\dots ,x_{n+1}]_d\subset \gamma :\
x\in \{x_1,\dots ,x_{n+1}\}}d_{\sigma ,x}^{X*}(W_x(\gamma ,x)_{[x_1,\dots
,x_{n+1}]_d}).  \label{5.5}
\end{equation}
Suppose, analogously to the above, that in local coordinates on the manifold
$X$
\begin{equation}
{\cal O}_{\gamma ,x}\ni y\mapsto W_x(\gamma ,y)_{[x_1,\dots
,x_{n+1}]_d}=w(y)\,h_1\wedge \dots \wedge h_{n+1},  \label{pils3}
\end{equation}
where $w\colon{\cal O}_{\gamma ,x}\to {\Bbb R}$ and $h_1\wedge \dots \wedge
h_{n+1}\in T_{x_1}X\wedge \dots \wedge T_{x_{n+1}}X$. Then,
\begin{multline}
d_{\sigma ,x}^{X*}(W_x(\gamma ,x)_{[x_1,\dots ,x_{n+1}]_d})=-\frac 1{\sqrt{%
n+1}}\sum_{i=1}^{n+1}(-1)^{i-1}\delta _{x,x_i}\big[ \langle \nabla
^Xw(x),h_i\rangle _x  \label{5.6} \\ \text{{}}+w(x)\langle \beta
_\sigma (x),h_i\rangle _x\big]h_1\wedge \dots \wedge
\check{h}{}_i\wedge \dots \wedge h_{n+1}.
\end{multline}
Here,
\[
\delta _{x,x_i}=%
\begin{cases}
1,&\text{if }x=x_i,\\
0,&\text{otherwise},\end{cases}
\]
and $\check{h}{}_i$ denotes the absence of $h_i$.

Upon \eqref{5.4}--\eqref{5.6}
\[
d_{\pi _\sigma }^{\Gamma *}\colon {\cal F}\Omega ^{n+1}\to
L^2_{\pi_\sigma}\Omega ^n.
\]

For $n\in{\Bbb N}$, we define the pre-Dirichlet form ${\cal E}%
^{\mathrm R}_{\pi_\sigma}$ by
\begin{multline}  \label{lklk}
{\cal E}^{\mathrm
R}_{\pi_\sigma}(W^{(1)},W^{(2)}):=\int_{\Gamma_X}\big[ \langle
d^\Gamma W^{(1)}(\gamma),d^\Gamma
W^{(2)}(\gamma)\rangle_{\wedge^{n+1}(T_\gamma\Gamma_X)} \\
\text{}+\langle
d^{\Gamma*}_{\pi_\sigma}W^{(1)}(\gamma),d^{\Gamma*}
_{\pi_\sigma}W^{(2)}(\gamma)\rangle_{\wedge^{n-1}(T_\gamma\Gamma_X)} \big]%
\,\pi_\sigma(d\gamma),
\end{multline}
where $W^{(1)},W^{(2)}\in {\cal F}\Omega^n$. Analogously to the
case of Bochner, we conclude that the function under the sign of
integral in \eqref{lklk} is polynomially bounded, so that the
integral exists.

The next theorem follows from \eqref{5.1}--\eqref{5.6}.

\begin{theorem}
\label{th5.1} For any $W^{(1)},W^{(2)}\in {\cal F}\Omega ^n${\rm , } we have
\[
{\cal E}_{\pi _\sigma }^{\mathrm R}(W^{(1)},W^{(2)})=\int_{\Gamma
_X}\langle H_{\pi _\sigma }^{\mathrm R}W^{(1)}(\gamma
),W^{(2)}\rangle _{\wedge ^n(T\Gamma _X)}\,\pi _\sigma (d\gamma ).
\]
Here{\rm , } $H_{\pi _\sigma }^ {\mathrm R}=d^\Gamma d_{\pi
_\sigma }^{\Gamma *}+d_{\pi _\sigma }^{\Gamma *}d$ is an operator
in the space $L_{\pi _\sigma }^2\Omega ^n$ with domain ${\cal
F}\Omega ^n${\rm . } It can be represented as follows{\rm : }
\begin{equation}
H_{\pi _\sigma }^{\mathrm R}W(\gamma )=\sum_{x\in \gamma
}H_{\sigma ,x}^{\mathrm R}W(\gamma )=\langle H_{\sigma ,\bullet
}^{\mathrm R}\,W(\gamma ),\gamma \rangle ,\qquad W\in {\cal
F}\Omega ^n, \label{5.7}
\end{equation}
where
\begin{equation}
H_{\sigma ,x}^{\mathrm R}=d_x^Xd_{\sigma ,x}^{X*}+d_{\sigma
,x}^{X*}d_x^X. \label{5.8}
\end{equation}
\end{theorem}

From Theorem~\ref{th5.1} we conclude that the pre-Dirichlet form ${\cal E}%
_{\pi _\sigma }^{\mathrm R}$ is closable in the space $L_{\pi
_\sigma }^2\Omega ^n$. The generator of its closure (being
actually the Friedrichs extension of the operator $H_{\pi _\sigma
}^{\mathrm R}$, for which we preserve the same notation) will be
called the de Rham Laplacian on $\Gamma _X$ corresponding to the
Poisson measure $\pi _\sigma $. By (\ref{5.7}) and (\ref{5.8}),
$H_{\pi
_\sigma }^{\mathrm R}$ is the lifting of the de Rham Laplacian on $X$ with measure $%
\sigma $.

\begin{remark}\label{remkll}
{\rm Similarly to (\ref{cyl-boch}), the operator $H_{\pi _\sigma
}^{\mathrm R}$ preserves the space ${\cal F}\Omega ^n$, and we can
always take $\Lambda (H_{\pi _\sigma }^{\mathrm R}W)=\Lambda (W)$.
Then, for any open bounded $\Lambda \supset \Lambda (W)$, we have
\begin{equation}
(H_{\pi _\sigma }^{\mathrm R}W)_{\Lambda ,\gamma }=H_{\sigma
^{\otimes |\Lambda \cap \gamma |}}^{\mathrm R}(X^{|\Lambda \cap
\gamma |})W_{\Lambda ,\gamma }, \label{cyl-der}
\end{equation}
where }$H_{\sigma ^{\otimes |\Lambda \cap \gamma |}}^B(X^{|\Lambda
\cap \gamma |})${\rm \ is the de Rham Laplacian of the manifold
$X^{|\Lambda \cap \gamma| }$ with the product measure $\sigma
^{\otimes|\Lambda \cap \gamma| }$.}
\end{remark}

Analogously to Theorem~\ref{thonsa1}, we get

\begin{theorem}
\label{thonsa2} {\rm 1)} On ${\cal D}\Omega ^n$ we have
\begin{equation}
H_{\pi _\sigma }^{\mathrm R}=(I^n)^{-1}\bigg[H_{\pi _\sigma }\boxplus \bigg( %
\bigoplus_{m=1}^nH_{\sigma ,\,(n,m)}^{\mathrm R}\bigg)\bigg]I^n,
\label{mn}
\end{equation}
where $H_{\sigma ,\,(n,m)}^{\mathrm R}$ denotes the restriction of
the de Rham Laplacian acting in the space $L^2(X^m\to \wedge
^n(TX^m);\sigma ^{\otimes m})$ to the subspace $L_\sigma ^2\Psi
^n(X^m)${\rm . }

{\rm 2)} Suppose that, for each $m=1,\dots ,n,$ the de Rham Laplacian $%
H_{\sigma ,m}^{\mathrm R}(X)$ is essentially self-adjoint on the
set of smooth forms with compact support. Then{\rm , } ${\cal
D}\Omega ^n$ is a domain of
essential self-adjointness of $H_{\pi _\sigma }^{\mathrm R}${\rm , } and the equality~%
\eqref{mn} holds for the closed operators $H_{\pi _\sigma
}^{\mathrm R}$ and $H_{\pi
_\sigma }\boxplus \big( \bigoplus_{m=1}^nH_{\sigma ,\,(n,m)}^{\mathrm R}\big)$ {\rm (}%
where the latter operator is closed from its domain of essential
self-adjointness $I^n({\cal D}\Omega ^n)${\rm )}{\rm . }
\end{theorem}

\begin{remark}\label{brand_new_remark2}\rom{ The essential
self-adjointness of the de Rham Laplacian $H^{\mathrm R}_{\sigma}$
on the set of smooth forms with compact support is well-known in
the case where $\sigma$ is the volume measure, see e.g.\
\cite{E3}. It is also  sufficient to assume that $\beta_\sigma$,
together with its derivatives up to order 3,  as well as the the
curvature tensor of $X$, together with its derivatives up to order
2, are bounded, cf.\ Remark~\ref{brand_new_remark1}. }\end{remark}

\noindent{\it Proof of Theorem}~\ref{thonsa2}. Upon \eqref{odyn}, \eqref{dvi}, \eqref{dunia}--%
\eqref{5.6}, \eqref{5.7}, and \eqref{5.8}, we get, for any $W\in{\cal D}%
\Omega^n$ given by the formula \eqref{kikimora},
\begin{gather*}
(H^{\mathrm R}_{\pi_\sigma}W)_k(\gamma)(\bar x)=0\qquad \text{for
}k\ne m, \\
(H^{\mathrm R}_{\sigma,x}W)_m(\gamma)(\bar x)=\begin{cases}(m!)^{1/2}\,(H_{\sigma,x}F)(%
\gamma\setminus\{\bar x\})\omega(\bar
x),&x\in\gamma\setminus\{\bar
x\},\\(m!)^{1/2}\,F(\gamma\setminus\{\bar x\})(H^{\mathrm
R}_{\sigma,x}\omega)(\bar x),&x\in\{\bar x\}.\end{cases}
\end{gather*}
Hence, analogously to \eqref{chyzh} and \eqref{pyzh}, we derive
\[
(I_k^n H^{\mathrm R}_{\pi_\sigma}W)(\gamma,\bar x)=%
\begin{cases}0,& k\ne m,\\
(H_{\pi_\sigma}F) (\gamma)\omega(\bar{x})+F(\gamma)(H^{\mathrm R
}_{\sigma,\,(n,m)}\omega)(\bar{x}),&k=m,\end{cases}
\]
which easily yields \eqref{mn}.

2)The proof is similar to that of Theorem~\ref{thonsa1}, 2). \quad $%
\blacksquare$

Again, analogously to Corollary~\ref{fock1}, we get a Fock space
representation of the operator $H^{\mathrm R}_{\pi_\sigma}$.

\begin{corollary}
\label{fock2} Let the conditions of Theorem~{\rm \ref{thonsa2}, 2)} be
satisfied{\rm . } Then{\rm , }
\begin{equation}
{\cal I}^nH_{\pi _\sigma }^{\mathrm R}({\cal
I}^n)^{-1}=d\operatorname{Exp}H_\sigma \boxplus \bigg(
\bigoplus_{m=1}^nH_{\sigma ,\,(n,m)}^{\mathrm R}\bigg).  \nonumber
\end{equation}
\end{corollary}

\subsection{Weitzenb\"{o}ck formula on the configuration  space}
\label{Subsec3.5}

In this section, we will derive a generalization of the Weitzenb\"{o}ck
formula to the case of the Poisson measure on the configuration space. In
other words, we will derive a formula which gives a relation between the
Bochner and de Rham Laplacians.

Analogously to \eqref{odyn}, \eqref{dvi}, we define for each $V(\gamma)\in
T_\gamma\Gamma_X$, $\gamma\in\Gamma_X$, the annihilation and creation
operators
\begin{align*}
a(V(\gamma))&\colon \wedge
^{n+1}(T_\gamma\Gamma_X)\to\wedge^n(T_\gamma\Gamma_X), \\
a^*(V(\gamma))&\colon \wedge^n
(T_\gamma\Gamma_X)\to\wedge^{n+1}(T_\gamma\Gamma_X)
\end{align*}
as follows:
\begin{alignat*}
{3} a(V(\gamma))W_{n+1}(\gamma)&=\sqrt{n+1}\,\langle
V(\gamma),W_{n+1}(\gamma)\rangle_\gamma,&&\qquad&&W_{n+1}(\gamma)
\in\wedge^{n+1}(T_\gamma\Gamma_X), \\
a^*(V(\gamma))W_n(\gamma)&=\sqrt{n+1}\,V(\gamma)\wedge W
_n(\gamma),&&\qquad&&W_n(\gamma)\in\wedge^n(T_\gamma\Gamma_X).
\end{alignat*}

Now, for a fixed $\gamma\in\Gamma_X$ and $x\in\gamma$, we define the
operator $R(\gamma)$ as follows:
\begin{gather*}
R(\gamma)=\sum_{x\in\gamma}R(\gamma,x),\qquad
D(R(\gamma)):=\wedge_0^n(T_\gamma\Gamma_X), \\
R(\gamma,x):=\sum_{i,j,k,l=1}^dR_{ijkl}(x)a^*(e_i)a(e_j)a^*(e_k)a(e_l).
\end{gather*}
Here, $\{e_j\}_{j=1}^d$ is a fixed orthonormal basis in the space $T_xX$
considered as a subspace of $T_\gamma\Gamma_X$, and $\wedge_0^n(T_\gamma%
\Gamma_X)$ consists of all $W(\gamma) \in\wedge^n(T_\gamma\Gamma_X)$ having
only a finite number of nonzero coordinates in the direct sum expansion %
\eqref{tang-n0}.

Next, we note that
\begin{align*}
\nabla^\Gamma B_{\pi_\sigma}(\gamma)&=(\nabla^X_x
B_{\pi_\sigma}(\gamma))_{x\in\gamma}=(\nabla^X_x(B_{\pi_\sigma}(\gamma)_y))
_{x,y\in\gamma} \\
&=(\delta_{x,y}\nabla^X\beta_\sigma(y))_{x,y\in\gamma}\in(T_{\gamma,\infty}%
\Gamma_X)^{\otimes2}.
\end{align*}
Hence, for any $V(\gamma)\in T_{\gamma,0}\Gamma_X$,
\begin{align*}
\nabla^\Gamma_V B_{\pi_\sigma}(\gamma):&= \langle\nabla^\Gamma
B_{\pi_\sigma}(\gamma),V(\gamma)\rangle_\gamma \\
&=\bigg( \sum_{y\in\gamma}\delta_{x,y}\langle\nabla^X\beta_\sigma(y),V(%
\gamma)_y\rangle_y \bigg)_{x\in\gamma} \\
&=\big(\langle\nabla^X\beta_\sigma(x),V(\gamma)_x\rangle_x\big)%
_{x\in\gamma}\in T_{\gamma,0}\Gamma_X .
\end{align*}
Thus, $\nabla^\Gamma B_{\pi_\sigma}(\gamma)$ determines the linear operator
in $T_{\gamma,0}\Gamma_X$ given by
\[
T_{\gamma,0}\Gamma_X\ni V(\gamma)\mapsto \nabla^\Gamma
B_{\pi_\sigma}(\gamma)V(\gamma):= \nabla^\Gamma_V B_{\pi_\sigma}(\gamma)\in
T_{\gamma,0}\Gamma_X.
\]
Analogously to \eqref{6.1}, we define in $\wedge_0^n(T_\gamma\Gamma_X)$ the
operator
\begin{multline*}
(\nabla^\Gamma B_{\pi_\sigma}(\gamma))^{\wedge n}:=\nabla^\Gamma
B_{\pi_\sigma} (\gamma)\otimes{\bf 1}\dots\otimes {\bf 1}+ {\bf 1}\otimes
\nabla^\Gamma B_{\pi_\sigma}(\gamma) \otimes {\bf 1}\otimes\dots\otimes{\bf 1%
} \\
\text{}+ \dots+{\bf 1}\otimes\dots\otimes{\bf 1}\otimes \nabla^\Gamma
B_{\pi_\sigma}(\gamma).
\end{multline*}

\begin{theorem}[Weitzenb\"ock formula on the Poisson space]
We have on ${\cal F}\Omega ^n$
\begin{equation}
H_{\pi _\sigma }^{\mathrm R}=H_{\pi _\sigma }^{\mathrm B}+R_{\pi
_\sigma }(\gamma ), \label{6.2}
\end{equation}
where
\begin{equation}
R_{\pi _\sigma }(\gamma ):=R(\gamma )-(\nabla ^\Gamma B_{\pi _\sigma
}(\gamma ))^{\wedge n}.  \label{6.2new}
\end{equation}
\end{theorem}

\noindent{\it Proof}. Fix $W\in{\cal F}\Omega^n$ and $\gamma\in\Gamma_X$.
Let $\Lambda(W)\subset X$ be a compactum as in Definition~\ref{def2.2}
corresponding to $W$, and let $\Lambda$ be an open set in $X$ with compact
closure such that $\Lambda(W)\subset\Lambda$. Next, let $W_{\Lambda,\gamma}$
be the form on ${\cal O}_{\gamma,x_1}\times\dots\times{\cal O}_ {\gamma,x_k}$%
, $\{x_1,\dots,x_k\}=\gamma\cap\Lambda$, defined by \eqref{cyl-form}.

It follows from Remarks~\ref{rem4.1}, 2) and \ref{remkll} that
\begin{align*}
\operatorname{Proj}_{\wedge ^n(T_{x_1}\oplus \dots \oplus
T_{x_k})}(H_{\pi _\sigma }^{\mathrm B}W(\gamma ))& =H_{\sigma
^{\otimes |\Lambda \cap \gamma |}}^{\mathrm B}(X^{|\Lambda \cap
\gamma |})W_{\Lambda ,\gamma }(x_1,\dots ,x_k), \\
\operatorname{Proj}_{\wedge ^n(T_{x_1}\oplus \dots \oplus
T_{x_k})}(H_{\pi _\sigma }^{\mathrm R}W(\gamma ))& =H_{\sigma
^{\otimes |\Lambda \cap \gamma |}}^{\mathrm R}(X^{|\Lambda \cap
\gamma |})W_{\Lambda ,\gamma }(x_1,\dots ,x_k),
\end{align*}
and $H_{\pi _\sigma }^{\mathrm B}W(\gamma )_{[y_1,\dots
,y_n]_d}=H_{\pi _\sigma }^{\mathrm R}W(\gamma )_{[y_1,\dots
,y_n]_d}=0$, $[y_1,\dots ,y_n]_d\subset \gamma $, if at least one
$y_i\in \{y_1,\dots ,y_n\}$ does not belong to $\Lambda $. Now,
the formulas \eqref{6.2}, \eqref{6.2new} follow from the usual
Weitzenb\"{o}ck formula \eqref{wei1}, \eqref{wei2} for the operators $%
H_{\sigma ^{\otimes |\Lambda \cap \gamma |}}^{\mathrm
B}(X^{|\Lambda \cap \gamma |})$ and $H_{\sigma ^{\otimes |\Lambda
\cap \gamma |}}^{\mathrm R}(X^{|\Lambda \cap \gamma |})$.\quad
$\blacksquare $\vspace{2mm}

We will show now that the Weitzenb\"{o}ck correction term
$R_{\pi_\sigma}$ is a lifting of the Weitzenb\"{o}ck correction
terms $R_{\sigma,k}$ of the manifold $X$.

Given  operator fields
\begin{equation}
X\ni x\mapsto J_k(x)\in {\cal L}(\wedge ^k(T_xX)),\qquad
k=1,\dots,\min\{n,d\},\label{chernysh}
\end{equation}
which are supposed to be uniformly bounded, we define a
``diagonal'' operator field
\begin{equation}
\widetilde X^m\ni \bar{x}\mapsto J_{n,m}(\bar{x})\in {\cal
L}({\Bbb T }_{\{\bar{x}\}}^{(n)}X^m),  \qquad
m=1,\dots,n,\label{field-j}
\end{equation}
as follows. First, we define for  each $\bar x=(x_1,\dots,x_m)\in
\widetilde X{}^m$ operators
\begin{gather*}
J_{n,m}^{k_1,\dots,k_m}(x_1,\dots,x_m)\in{\cal
L}((T_{x_1}X)^{\wedge k_1}\wedge\dots\wedge(T_{x_m}X)^{k_m}),\\
1\le k_1,\dots,k_m\le d,\quad k_1+\dots+k_m=n,
\end{gather*}
by setting
\begin{gather}
J_{n,m}^{k_1,\dots,k_m}(x_1,\dots,x_m)
u^{(k_1)}_1\wedge\dots\wedge u^{(k_m)}_m:=
\big(J_{k_1}(x_1)u^{(k_1)}_1\big)\wedge
u^{(k_2)}_2\wedge\dots\wedge u^{(k_m)}_m \notag\\
\text{}+u^{(k_1)}_1\wedge\big(J_{k_2}(x_2)u^{(
k_2)}_2\big)\wedge\dots\wedge u^{(k_m)}_m
+\dots+u^{(k_1)}_1\wedge\dots\wedge u^{(k_{m-1})}_ {m-1}
\wedge\big(J_{k_m}(x_m)u^{(k_m)}_m\big),\notag\\
u^{(k_i)}_i\in\wedge^{k_i}(T_{x_i}X),\
i=1,\dots,m.\label{field-j1}
\end{gather}
and  extending the operator
$J_{n,m}^{k_1,\dots,k_m}(x_1,\dots,x_m)$ by linearity and
continuity to the whole space. Then, the operator
$J_{n,m}(x_1,\dots,x_m)\in{\cal L}({\Bbb
T}^{(n)}_{\{x_1,\dots,x_m\}}X^m)$ is defined by setting its
diagonal blocks in the decomposition (\ref {tang-n2}) of the space
${\Bbb T}^{(n)}_{\{x_1,\dots,x_m\}}X^m$ to be
$J_n^{k_1,\dots,k_m}(x_1,\dots,x_m)$ and the other blocks to be
equal to zero.

Notice that, for each $\nu\in S_m$, the operators
$J_{n,m}(x_1,\dots,x_m)$ and
$J_{n,m}(x_{\nu(1)},\dots,x_{\nu(m)})$ coincide, so that
\eqref{field-j} naturally determines the operator field
\begin{equation}
\widetilde X^m/S_m\ni \{\bar{x}\}\mapsto J_{n,m}(\{\bar{x}\})\in
{\cal L}({\Bbb T }_{\{\bar{x}\}}^{(n)}X^m),  \qquad
m=1,\dots,n,\label{hdfushdf}
\end{equation}

Now, we define an operator field
\begin{equation}
{\bf \Gamma }_X\ni \gamma \mapsto {\bf J}(\gamma )\in {\cal
L}(\wedge ^n(T_\gamma \Gamma _X))
\end{equation}
setting ${\bf J}(\gamma )$ to be again the block-diagonal operator
in the decomposition (\ref{tang-n1}) with the diagonal blocks
$J_{n,m}(\{\bar{x}\})$ and the other blocks equal to zero.

In what follows, we suppose, for simplicity, that
\begin{equation}\label{ninzia}
\text{the curvature tensor $R_{ijkl}(x)$ and
$\nabla^X\beta_\sigma(x)$ are uniformly bounded in $x\in
X$.}\end{equation} As easily seen, for each $k\in\N$, the
Weitzenb\"ock correction term $R_{\sigma,k}(\bullet)$ on the
manifold $X$ is now a uniformly bounded operator field taking
values in  $\wedge^k(TX)$ (cf. \ref{wei2}). Thus we can define an
operator field ${\bf R}_\sigma (\gamma )$ through  the operator
fileds $ R_{\sigma,k} (x)$.

\begin{proposition} Let \eqref{ninzia} hold\rom.
Then\rom,
\[
R_{\pi _\sigma }={\bf R}_\sigma .
\]
\end{proposition}

\noindent {\it Proof}. The result can be easily seen directly from the definition of $%
{\bf R}_\sigma(\gamma)$, $R(\gamma)$ and $B_{\pi _\sigma }(\gamma
)$.\quad $\blacksquare $

\section{Probabilistic representation of the Laplacians}

Let $\xi _x(t)$ be the Brownian motion on $X$ with the drift $\beta _\sigma $---
the logarithmic derivative of $\sigma $---which starts at a point
$x\in X$. We suppose the following:

\begin{itemize}
\item  for each $x\in X$, the process $\xi _x(t)$ has  infinite life-time;

\item  the semigroup
\begin{equation}\notag
T_0(t)f(x)%
:=%
{\sf E}\,f(\xi _x(t))
\end{equation}
preserves the space $C_{\mathrm b}^2(X)$ and can be extended to a
strongly continuous semigroup of contractions in $L^2(X;\sigma )$,
and its generator $H_0$ is essentially self-adjoint on the space
$\cal D$  (in this case $H_0=H_\sigma $).
\end{itemize}

\begin{remark}\rom{
The above conditions are fulfilled if e.g.\ $\beta _\sigma $,
together with its derivatives up to order $3$, is bounded.}
\end{remark}

We denote by $\xi _\gamma (t)$ the corresponding independent
infinite particle
process which starts at a point $\gamma \in \Gamma _X$,%
\begin{equation}\notag
\xi _\gamma (t)=(\xi _x(t))_{x\in \gamma }.
\end{equation}
As we already mentioned in subsec.~\ref{subsec3.1}, this process
is properly associated with the Dirichlet form ${\cal E}_{\pi
_\sigma } $, see \cite{AKR1}.

\begin{remark}\rom{
The process $\xi _\gamma $  lives in general on the bigger state
space $\ddGamma$ consisting of all ${\Bbb Z}_{+}$-valued Radon
measures on $X$ (the space  $\ddGamma$ being Polish). Notice,
however, that at each fixed moment of time $t\in{\Bbb R}_+$ the
value $\xi _\gamma (t)$ belongs to $\Gamma _X$ a.s. Moreover, it
was proven in \cite{RSm} that, in the special case $X={\Bbb R}^d$
with $d\ge 2$, the process $\xi _\gamma $ lives a.s.\ in $\Gamma
_X$.}
\end{remark}

Let ${\bf T}_0(t)F(\gamma )%
:=%
{\sf E}\,F(\xi _\gamma (t))$ be the corresponding semigroup. As
shown in \cite{AKR1}, it can be extended from ${\cal FC}_{\mathrm
b}^\infty
(\Gamma _X)$ to a strongly continuous semigroup in $%
L_{\pi _\sigma }^2(\Gamma _X)$ with the generator ${\bf
H}_0=H_{\pi _\sigma }^\Gamma $.

Given  operator fields \eqref{chernysh}
 which are now supposed to be uniformly bounded,
continuous,  and symmetric, we define again  operator fields
(\ref{field-j}) in the same way as in subsec.~\ref{Subsec3.5}.
 We have obviously
$J_{n,m}(\bar{x})^{*}=J_{n,m}(\bar{x})$.

Let
\begin{equation}\notag
P_{\xi _{\bar{x}}}^{J_{n,m}}(t)\colon{\Bbb
T}_{\{\bar{x}\}}^{(n)}X^m\to {\Bbb T}_{\{\xi
_{\bar{x}}(t)\}}^{(n)}X^m,\qquad m=1,\dots,n
\end{equation}
be the parallel translation along the path $\xi _{\bar{x}}(t)%
:=%
(\xi _{x_i}(t))_{i=1,\dots,m}$ with the potential $J_{n,m}$. That
is, $\eta (t)=P_{\xi _{\bar{x}}}^{J_{n,m}}(t)h$ satisfies the SDE
\begin{equation}
\frac D{dt}\,\eta (t)=J_{n,m}(\eta (t)),\qquad\eta (0)=h,
\end{equation}
where $\frac D{dt}$ denotes the covariant differentiation along
the paths of the process $\xi $ (see \cite{E3}). It is easy to see
that the symmetry of the potential  $J_{n,m}(\bar x)$ w.r.t.\ a
permutation of the components of $\bar x$ implies the same
symmetry of $P_{\xi _{\bar{x}}}^{J_{n,m}}(t)$. Thus, analogously
to \eqref{hdfushdf}, we get the operator field
\begin{equation}\label{ghjklo}
\widetilde{X}^m/S_m\ni \{\bar{x}\}\mapsto P_{\{\xi _{\bar{x}%
}\}}^{J_{n,m}}(t). \end{equation}

Now, for $\pi_\sigma$-a.e.\ $\gamma \in \Gamma _X$, we define the
operator
\begin{equation}\notag
{\bf P}_{\xi _\gamma }^{\bf J}(t)\colon\wedge ^n(T_\gamma \Gamma
_X)\to\wedge ^n(T_{\xi _\gamma (t)}\Gamma _X)
\end{equation}
by setting its diagonal blocks in the decomposition (\ref{tang-n1}) to be $%
P_{\{\xi _{\bar{x}}\}}^{J_{n,m}}(t)$ and the other blocks to be
equal to zero.

It is known that
\begin{equation}
\| P_{\xi _{\bar{x}}}^{J_{n,m}}(t)\| \le e^{tC_m},\qquad
m=1,\dots,n,  \label{bound1}
\end{equation}
where $C_m$ is the supremum of the spectrum of $J_{n,m}(\bar{x})$.

\begin{lemma}
For $\pi_\sigma$-a\rom.e\rom.\ $\gamma \in \Gamma _X$\rom, we have
\begin{equation}
\| {\bf P}_{\xi _\gamma }^{\bf J}(t)\| \le e^{tC},\qquad
C=\max_{m=1,\dots,n}C_m. \label{est1}
\end{equation}
\end{lemma}

\noindent {\it Proof}. The result follows directly from the definition of ${\bf P}%
_{\xi _\gamma }^{\bf J}(t)$ and estimate (\ref{bound1}).\quad $
\blacksquare $\vspace{2mm}

Let us define a semigroup ${\bf T}_n^{{\bf J}}(t)$ acting in the space of $%
n$-forms as follows:
\begin{equation}
{\bf T}_n^{{\bf J}}(t)W(\gamma )%
:=%
{\sf E}\,\big( {\bf P}_{\xi _\gamma }^{\bf J}(t)\big) ^{*}W(\xi
_\gamma (t)),\qquad W\in{\cal F}\Omega^n.
\end{equation}

Let $T_{n,m}^J(t)$ be the semigroup acting in the space $L_\sigma
^2\Psi ^n(X^m)$ as
\begin{equation}
T_{n,m}^J(t)\omega (\bar{x})%
:=%
{\sf E}\,\big(P_{\xi _{\bar{x}}}^{J_{n,m}}(t)\big)^{*}\omega (\xi _{\bar{x}%
}(t)).
\end{equation}
By virtue of \eqref{ghjklo} and estimate (\ref{bound1}), we
conclude  the correctness of the definition of $T_{n,m}^J(t)$ (in
the sense that  $T_{n,m}^J(t)$ is uniquely defined) and its strong
continuity. The following
result describes the structure and properties of the semigroup ${\bf T}_n^{%
{\bf J}}(t)$.

\begin{proposition}
\label{pnsem} \rom{1)}  ${\bf T}_n^{{\bf J}}(t)$ satisfies the
estimate
\begin{equation}
\| {\bf T}_n^{{\bf J}}(t)V(\gamma )\| _{\wedge ^n(T_\gamma \Gamma
_X)}\le e^{tC}\,{\bf T}_0(t)\| V(\gamma )\| _{\wedge ^n(T_\gamma
\Gamma _X)}  \label{markov}
\end{equation}
for $\pi_\sigma$-a\rom.e\rom.\ $\gamma \in \Gamma _X.$

\rom{2)}  Under the  isomorphism $I^n$\rom, ${\bf T}_n^{{\bf
J}}(t) $ takes the following form\rom:
\begin{equation}
I_m^n{\bf T}_n^{{\bf J}}(t)={\bf T}_0(t){\bf \otimes
}T_{n,m}^J(t)\;I_m^n, \qquad m=1,\dots,n.\label{dec-sem}
\end{equation}
In particular\rom, for $1$-forms
\begin{equation}
I^1{\bf T}_1^{{\bf J}}(t)={\bf T}_0(t){\bf \otimes
}T_{1,1}^J(t)I^1.
\end{equation}

\rom{3)} ${\bf T}_n^{{\bf J}}(t)$ extends to a strongly continuous
semigroup in $L_{\pi _\sigma }^2\Omega ^n.$

\end{proposition}

\noindent {\it Proof}. 1)  The result follows from formula
(\ref{est1}).

2)  For simplicity, we give the proof only in the case of
$1$-forms. Let $V\in {\cal D}\Omega ^1$ be given by $I^1V=F\otimes
v$. By the
definition of ${\bf T}_1^{{\bf J}}(t)$ and the construction of the process $%
\xi _\gamma $, we have
\begin{align*}
{\bf T}_1^{{\bf J}}(t)V(\gamma )_x &={\sf E\,}F(\xi _\gamma
(t)\setminus \left\{ \xi _x(t)\right\} )\big(P_{\xi
_x}^{J_{1,1}}(t)\big)^{*}v(\xi _x(t))
\\
&={\sf E\,}F(\xi _\gamma (t)\setminus \left\{ \xi _x(t)\right\} ){\sf E}%
_{\xi _x}\,\big(P_{\xi _x}^{J_{1,1}}(t)\big)^{*}v(\xi _x(t))
\\  &={\bf T}_0(t)F(\gamma \setminus \left\{
x\right\} )T_1^J(t)v(x),
\end{align*}
${\sf E}_{\xi _x}$ meaning the expectation w.r.t.\ the process
$\xi _x(t)$, from where the result follows. The general case can
be proved by similar arguments.

3)  The result follows from the corresponding results for semigroups $%
{\bf T}_0(t)$ and $T_{n,m}^J(t)$, which are well-known (see \cite
{AKR-1}, resp. \cite{E3}).\quad $ \blacksquare $\vspace{2mm}

Let ${\bf H}_n^{{\bf J}}$ and $H_{n,m}^J$ be the generators of ${\bf T}_n^{%
{\bf J}}(t)$ and $T_{n,m}^J(t)$, respectively.

Now, we will give probabilistic representations of the semigroups
$T_{\pi _\sigma }^{\mathrm B}(t)$ and $T_{\pi _\sigma }^{\mathrm R
}(t)$ associated with the operators $H_{\pi _\sigma }^{\mathrm B}$
and $H_{\pi _\sigma }^{\mathrm R}$, respectively. We set
$$J_m^{(1)}:=0,\qquad J_m^{(2)}(x):=R_{\sigma ,m}(x),\qquad
m=1,\dots,\min \{n,d\}$$ (cf.\ (\ref{wei2})). Let us remark that
$P_{\xi _{\bar{x}}}^{J_{n,m}^{(1)}}(t)\equiv P_{\xi
_{\bar{x}}}(t)$ is the parallel translation of the $n$-forms along the path $%
\xi _{\bar{x}}$, and we have $$H_{n,m}^{J^{(1)}}=-H_{\sigma
,\,(n,m)}^{\mathrm B},\qquad H_{n,m}^{J^{(2)}}=-H_{\sigma
,\,(n,m)}^{\mathrm R}\qquad\text{on }\Psi _0^n(X^m).$$

\begin{theorem}
{\rm 1)} For $W\in {\cal D}\Omega ^n${\rm ,}  we have
\begin{equation}
H_{\pi _\sigma }^{\mathrm B}W=-{\bf H}_n^{{\bf J}^{(1)}}W,\qquad H_{\pi _\sigma }^
{ \mathrm R}W=-%
{\bf H}_n^{{\bf J}^{(2)}}W.  \label{prrepgen}
\end{equation}

{\rm 2)} As $L^2$-semigroups{\rm , }
\begin{equation}
T_{\pi _\sigma }^{\mathrm B}(t)={\bf T}_n^{{\bf
J}^{(1)}}(t),\qquad T_{\pi _\sigma }^{\mathrm R}(t)={\bf
T}_n^{{\bf J}^{(2)}}(t). \label{prrepsem}
\end{equation}

{\rm 3)} The semigroups $T_{\pi _\sigma }^{\mathrm B}(t)$ and
$T_{\pi _\sigma }^{\mathrm R}(t)$ satisfy the estimates
\[
\Vert T_{\pi _\sigma }^{\mathrm B}(t)V(\gamma )\Vert _\gamma \le
{\bf T}_0(t)\Vert V(\gamma )\Vert _\gamma
\]
and
\[
\Vert T_{\pi _\sigma }^{\mathrm R}(t)V(\gamma )\Vert _\gamma \le e^{tC}\,{\bf T}%
_0(t)\Vert V(\gamma )\Vert _\gamma
\]
for $\pi_\sigma$-a\rom.e\rom.\ $\gamma \in \Gamma _X$\rom.
\end{theorem}

\noindent {\it Proof}. 1) It follows directly from the
decomposition (\ref {dec-sem}) that, on ${\cal D}\Omega ^n$, we
have
\begin{equation}
I_m^n\,{\bf H}_n^{{\bf J}}=\left( {\bf H}_0\boxplus
H_{n,m}^J\right) \,I_m^n ,\label{dec-gen}
\end{equation}
where ${\bf H}_0$ is the generator of ${\bf T}_0(t)$. Setting
respectively $J_m:=J_m^{(1)}$ and $J_m:=J_m^{(2)}$ and comparing
(\ref{dec-gen0}) with (\ref{dec-gen}), we obtain the result.

2) The statement follows from (\ref{prrepgen}) and the essential
self-adjointness of $H_{\pi _\sigma }^{\mathrm B}$ and $H_{\pi _\sigma }^{\mathrm R}$
on ${\cal D%
}\Omega ^n$ by applying Proposition~\ref{pnsem}, 3), with
$J_m=J_m^{(1)}$ and $J_m=J_m^{(2)}$.

3) The result follows from (\ref{prrepsem}) and (\ref{markov}).\quad $%
\blacksquare $

\section*{Acknowledgments}

It is a great pleasure to thank our friends and colleagues Yuri
Kondratiev, Tobias Kuna, and Michael R\"{o}ckner for their
interest in this work and the joy of collaboration on related
subjects. We are also grateful to V. Liebscher for a useful
discussion. The financial support of SFB 256, DFG Research Project
AL 214/9-3, and BMBF Project~UKR-004-99 is gratefully
acknowledged.

\end{document}